\documentclass[a4paper,twoside,11pt,reqno]{amsart}
\newcommand{\Ueberschrift}{
Weil--Ch\^atelet divisible elements in Tate--Shafarevich groups I: \\ The Bashmakov problem for elliptic curves over $\bQ$} 
\newcommand{\Kurztitel}{The Bashmakov problem for elliptic curves over $\bQ$}
\usepackage[headings]{fullpage}
\usepackage{amsmath} \usepackage{amssymb} \usepackage{amsopn}
\usepackage{amsthm} \usepackage{amsfonts} \usepackage{mathrsfs} 
\usepackage[dvips]{graphicx} 
\usepackage[all]{xy}
\usepackage[OT2, T1]{fontenc}
\usepackage[pdfpagelabels, pdftex]{hyperref}
\hypersetup{
  pdftitle={},
  pdfauthor={},
  pdfsubject={},
  pdfkeywords={},
  colorlinks=true,
  breaklinks=true,
  bookmarksopen=true,
  bookmarksnumbered=true,
  pdfpagemode=UseOutlines,
  plainpages=false}

\DeclareMathOperator{\rH}{H}
\DeclareMathOperator{\rI}{I}

\DeclareMathOperator{\rN}{N}

\DeclareMathOperator{\rc}{c}

\newcommand{\bF}{{\mathbb F}}

\newcommand{\bL}{{\mathbb L}}

\newcommand{\bN}{{\mathbb N}}

\newcommand{\bP}{{\mathbb P}}
\newcommand{\bQ}{{\mathbb Q}}

\newcommand{\bT}{{\mathbb T}}

\newcommand{\bZ}{{\mathbb Z}}

\newcommand{\cF}{{\mathscr F}}

\newcommand{\cM}{{\mathscr M}}

\newcommand{\fm}{{\mathfrak m}}

\newcommand{\fo}{{\mathfrak o}}
\newcommand{\fp}{{\mathfrak p}}

\DeclareSymbolFont{cyrletters}{OT2}{wncyr}{m}{n}
\DeclareMathSymbol{\Sha}{\mathalpha}{cyrletters}{"58}

\newcommand{\one}{\mathbf{1}}
\newcommand{\surj}{\twoheadrightarrow} 
\newcommand{\inj}{\hookrightarrow}

\newcommand{\ev}{{\rm ev}}
\DeclareMathOperator{\Hom}{Hom}
\DeclareMathOperator{\Aut}{Aut}
\DeclareMathOperator{\End}{End}

\DeclareMathOperator{\im}{im}

\DeclareMathOperator{\GL}{GL}
\DeclareMathOperator{\PGL}{PGL}
\DeclareMathOperator{\SL}{SL}

\DeclareMathOperator{\coRk}{corank}

\newcommand{\matzz}[4]{\left(
\begin{array}{cc} #1 & #2 \\ #3 & #4 \end{array} \right)}

\DeclareMathOperator{\Spec}{Spec}
\DeclareMathOperator{\Div}{Div}
\DeclareMathOperator{\divisor}{div}
\DeclareMathOperator{\Pic}{Pic}

\DeclareMathOperator{\Frob}{Frob}

\DeclareMathOperator{\res}{res}
\DeclareMathOperator{\infl}{inf}
\DeclareMathOperator{\Gal}{Gal}

\newcommand{\ph}{\varphi}
\newcommand{\et}{\text{\rm \'et}} 
\newcommand{\alg}{{\rm alg}}
\newcommand{\nr}{{\rm nr}} 
\newcommand{\ab}{{\rm ab}} 
\newcommand{\Sel}{{\rm Sel}} 
\newcommand{\kum}{{\rm kum}}

\newcommand{\bruch}[2]{\genfrac{}{}{0.5pt}{}{#1}{#2}}
\newcommand{\ov}[1]{\mbox{${\overline{#1}}$}} 

\newtheorem{thm}{Theorem}
\newtheorem{prop}[thm]{Proposition}
\newtheorem{lem}[thm]{Lemma}
\newtheorem{cor}[thm]{Corollary}

\newtheorem*{thmA}{Theorem A}
\newtheorem*{thmB}{Theorem B}

\theoremstyle{remark}
\newtheorem{rmk}[thm]{Remark}
\newtheorem{ex}[thm]{Example}

\newenvironment{pro*}[1][Proof]{{\it{#1:}} }{}

\renewcommand{\labelenumi}{(\arabic{enumi})}

\numberwithin{equation}{section}

\begin{document}

\title[\Kurztitel]{\Ueberschrift} 
\author{Mirela \c{C}iperiani}
\address{Mirela \c{C}iperiani, Department of Mathematics, the University of Texas at Austin, 1 University Station, C1200 
Austin, Texas 78712, USA}
\email{mirela@math.utexas.edu}
\urladdr{http://www.ma.utexas.edu/users/mirela/}

\author{Jakob Stix}
\address{Jakob Stix, Mathematisches  Institut, Universit\"at Heidelberg, Im Neuenheimer Feld 288, 69120 Heidelberg, Germany}
\email{stix@mathi.uni-heidelberg.de}
\urladdr{http://www.mathi.uni-heidelberg.de/~stix/}
\thanks{The authors acknowledge the hospitality and support provided by MATCH and the Newton Institute. \\
\indent The first author was partially supported by an NSF and an NSA grant during the preparation of this manuscript.}
\date{September 5, 2012} 

\maketitle

\begin{quotation} 
  \noindent \small {\bf Abstract} --- For an abelian variety $A$ over a number field $k$ we discuss the maximal divisibile subgroup of $\rH^1(k,A)$ and its intersection with the subgroup $\Sha(A/k)$. The results are most complete for elliptic curves over $\bQ$. 
\end{quotation}


\setcounter{tocdepth}{1} {\scriptsize \tableofcontents}


\section{Introduction}


\subsection{The Bashmakov problem} 
Let $A/k$ be an abelian variety over an algebraic number field $k$ with algebraic closure $k^\alg$. Bashmakov  \cite{bashmakov:div,bashmakov:survey}, studied the question of whether elements of the Tate--Shafarevich group 
\[
\Sha(A/k)
\]
become divisible in the  Weil--Ch\^atelet group $\rH^1(k,A)= \rH^1(k,A(k^\alg))$, i.e., lie in the subgroup 
\[
\divisor(\rH^1(k,A))
\]
of divisible elements.  This question was initially asked by Cassels in the case of elliptic curves (see \cite{cassels:III}  Problem 1.3) because an affirmative answer would prove that the kernel of the Cassels' pairing equals the maximal divisible subgroup of the Tate-Shafarevich group.  

Bashmakov \cite{bashmakov:div} also investigates whether $\Sha(A/k)$ can meet the maximal divisible subgroup of $\rH^1(k,A)$,
\[
\Div(\rH^1(k,A)),
\]
 in a nontrivial way, see also \cite{harariszamuely:galsec}  \S4. Bashmakov's results are recalled in Section \S\ref{sec:bashmakov}. 
In other words, the Cassels--Bashmakov problem considers  the filtration
\begin{equation} \label{eq:filtration}
\rH^1(k,A)  \supseteq \divisor(\rH^1(k,A)) \supseteq \Div(\rH^1(k,A)) \supseteq 0
\end{equation}
and intersects it with $\Sha(A/k)$.  Observe that $\divisor(\rH^1(k,A))$ is generally strictly larger than the maximal divisible subgroup $\Div(\rH^1(k,A))$, see Section \S\ref{def-divs} and Section \S\ref{sec:selmerexample}.

\smallskip

In this article we focus on Bashmakov's question on the intersection of the Tate-Shafarevich group and the maximal divisible subgroup of the Weil-Ch\^atelet group. We address Cassels' original question in a separate article \cite{ciperianistix:cassels}. 

\smallskip

Our motivation for studying the Cassels--Bashmakov problem is twofold. The group  $\Sha(A/k)$ conjecturally has no divisible elements,  and this is known for elliptic curves over $\bQ$ of analytic rank $\leq 1$. The study of $\Sha(A/k)$ in the bigger group $\rH^1(k,A)$ may shed some light on the structure of $\Sha(A/k)$ itself in the general case. 

Secondly, the Bashmakov problem arises naturally in anabelian geometry when one wants to get hold on the index of a hyperbolic curve and therefore passes to the abelianization in form of the universal Albanese torsor. The connection\footnote{The authors acknowledge the influence of work of Harari and Szamuely for bringing the results of Bashmakov to their attention and for informing the first author about their value from the point of view of anabelian geometry.} more precisely is as follows. Let $W$ be a principal homogeneous space of $A$ over $k$, with a geometric point $\bar{w}$ of the base change $\ov{W} = W \times_k k^\alg$. Let $\Gal_k = \Gal(k^\alg/k)$ be the absolute Galois group $k$. A section of the fundamental exact sequence
\[
1 \to \pi_1(\ov{W},\bar{w}) \to \pi_1(W,\bar{w}) \to \Gal_k \to 1,
\]
yields that the corresponding class $[W] \in \rH^1(k,A)$ lies in the maximal divisible subgroup, see \cite{harariszamuely:galsec} and also \cite{stix:habil} Remark 176(3). If, in addition, local points exist on $W$, then $[W]$ belongs to $\Sha(A/k)$ and the existence of a global $k$-rational point follows from a negative answer to the Bashmakov problem for $A/k$, see again \cite{harariszamuely:galsec}  \S4 and also \cite{stix:habil} \S13.


\subsection{Summary of results} 
In view of the conjectured finiteness of $\Sha(A/k)$ the triviality of the $p$-primary part of
\[
\Div(\rH^1(k,A)) \cap \Sha(A/k).
\]
should be guaranteed for large primes $p$ depending on $A/k$.  Our aim therefore is to identify conditions on  a prime number $p$ which imply the triviality of the  
$p$-primary part of the intersection of the Tate-Shafarevich group with the maximal divisible subgroup of the Weil-Ch\^atelet group without relying on the triviality of 
$\Sha(A/k)_p$.

When we can answer the Bashmakov problem in the negative, i.e., when we can prove that $\Div(\rH^1(k,A))$ intersects $\Sha(A/k)$ trivially, our method actually shows more. In order to present the results, we define the \textbf{locally divisible $\rH^1$} as the kernel
\[
\rH^1_{\divisor}(k,A) = \ker\Big(\rH^1(k,A) \to \bigoplus_v \rH^1(k_v,A)/\divisor\big(\rH^1(k_v,A)\big) \Big).
\]
Since $\divisor(-)$ is a functor, we find 
\[
\divisor(\rH^1(k,A)) \subseteq \rH^1_{\divisor}(k,A).
\]
Focusing on the $p$-part and using local Tate duality
$\rH^1(k_v,A) = \Hom(A^t(k_v),\bQ/\bZ)$
we find for $v \nmid p$ that $\divisor(\rH^1(k_v,A)_{p^\infty}) = 0$ and therefore an exact sequence
\begin{equation} \label{eq:pparth1localdiv}
0 \to \Sha(A/k)_{p^\infty} \to \rH^1_{\divisor}(k,A)_{p^\infty} \to \bigoplus_{v \mid p} \divisor\big(\rH^1(k_v,A)_{p^\infty}\big) \cong (\bQ_p/\bZ_p)^{\dim(A) \cdot [k:\bQ]}.
\end{equation}

Concerning the relative position of $\Sha(A/k)$ and $\Div(\rH^1(k,A))$ our main result is for an elliptic curve $E/\bQ$ under the assumption of trivial analytic rank (see Sect.~\S\ref{sec:maxdivEllQ}). The following theorem is a summary of the results in this direction 
(part (4) follows essentially from  \cite{harariszamuely:galsec} Cor 4.2.).

\begin{thmA}
Let $E/\bQ$ be an elliptic curve.
\begin{enumerate}
\item If $E/\bQ$ has trivial algebraic rank, then $\Div(\rH^1(\bQ,E))$ contains a copy of $\bQ/\bZ$.
\item If $E/\bQ$ has trivial analytic rank, then $\Div(\rH^1(\bQ,E)) = \bQ/\bZ$ and, for odd primes $p$ of good reduction such that the $\Gal_\bQ$-representation on the $p$-torsion $E_p$ is irreducible, we find
\[
\Sha(E/\bQ)_{p^\infty} \oplus \Div\big(\rH^1(\bQ,E)\big)_{p^\infty} = \rH^1_{\divisor}(\bQ,E)_{p^\infty}.
\]
\item If $E/\bQ$ has positive algebraic rank, then $\Div(\rH^1(\bQ,E)) = \Div(\Sha(E/\bQ))$.
\item If $E/\bQ$ has analytic rank $1$, then $\Div(\rH^1(\bQ,E)) = 0$.
\end{enumerate}
\end{thmA}

\begin{rmk}
Observe that the $\Gal_\bQ$-representation on $E_p$ is irreducible for all
\[
p\neq 2,3,5,7,11,13,17, 19, 37, 43, 67, 163
\]
and for all $p>7$ if $E$ is semistable (see \cite{mazur:eisenstein}).
\end{rmk}

With respect to the filtration \eqref{eq:filtration} we deduce that even if  $\Div(\rH^1(\bQ,E))$ is non-trivial, it  can be strictly smaller than $\divisor(\rH^1(\bQ,E))$.
Moreover, it turns out that an old example of Selmer provides an example where the Tate--Shafarevich group intersects the maximal divisible group nontrivially, see Section~\S\ref{sec:selmerexample} for the following.

\begin{thmB} The Tate-Shafarevich group $\Sha(E/\bQ)$ can intersect  $\Div(\rH^1(\bQ,E))$ non-trivially for an elliptic curve $E/\bQ$ of trivial analytic rank. In particular, the Selmer curve  
\[
3X^3+4Y^3+5Z^3 =0
\]
represents a nontrivial class in 
\[
\Sha(E/\bQ) \cap \Div(\rH^1(\bQ,E))
\] 
with $E$ the Jacobian of the Selmer curve.
\end{thmB}

\subsection{Plan of the paper} 

In Section~\S\ref{sec:prelim} and Section~\S\ref{sec:ec} we discuss background from Galois cohomology and \'etale cohomology. 
The material of Section~\S\ref{sec:ec} is only used in Section~\S\ref{sec:bashmakov}, where we discuss Bashmakov's work, and in Section~\S\ref{sec:selmerexample}, where we prove Theorem B. Theorem A is proven in Section~\S\ref{sec:maxdiv} based on  the computation of a generalized Selmer group done in Section~\S\ref{sec:induction} which in turn depends on some group theory in $\GL_2(\bF_\ell)$ that is discussed in Section~\ref{sec:group}.


\smallskip

\subsection*{Acknowledgments} 
The authors would like to thank Dorian Goldfeld, David Harari, Ken Ribet and Tam\'as Szamuely for several useful discussions. The first author is also grateful to her advisor, Andrew Wiles, for introducing her to this method of thinking about the $p$-divisibility of the Tate-Shafarevich group. 
Finally, we would also like to thank the referees for several helpful suggestions.

\smallskip


\subsection{Notation}
We fix some notation which will be in use throughout the text.

\subsubsection{On number fields} \label{sec:curves}
In the sequel, we let $k$ be an algebraic number field, i.e., a finite extension of $\bQ$. 
The completion of $k$ at a place $v$ is $k_v$.
For a finite $\Gal_k$-module $M$ the field extension $k(M)/k$ is the fixed field of $k^\alg$ under $\ker(\Gal_k \to \Aut(M))$.

We set $X = \Spec(\fo_k)$ with the ring of integers $\fo_k$ in $k$. 
The finite places of $k$ will be identified with the closed points of $X$. 
An abelian variety $A/k$ with good reduction over a nonempty open $U \subset X$ extends to an abelian scheme over $U$ that we denote by $A/U$ by abuse of notation.

\subsubsection{On abelian groups} \label{def-divs}Let $M$ be an abelian group. The $n$-torsion subgroup is denoted by by $M_n$, and $M_{p^\infty}$ denotes the $p$-primary torsion $\bigcup_n M_{p^n}$. The subgroup 
\[
\divisor(M)=\bigcap_{n\geq 1} nM
\]
of divisible elements of $M$ contains the maximal divisible subgroup 
\[
\Div(M)
\]
of $M$, which equals $\bigcup \im(\ph)$ for all $\ph : \bQ \to M$. For matters of clarity we will distinguish between \textbf{"$a$ is divisible by $p$"}, meaning there is $a'$ with $a = p a'$, and \textbf{"$a$ is $p$-divisible"}, which means that for every $n \geq 1$ there is $a'$ with $a = p^n a'$.
For a discussion of $\divisor(M)$ and $\Div(M)$ see Jannsen \cite{jannsen:continuous} \S4. In particular note that in general $\divisor(M)$ is strictly larger than $\Div(M)$, as in the following example: 
\[
M= \big(\bigoplus_n \frac{1}{2n}\bZ/\bZ\big)/\ker \big(\text{sum}:\bigoplus_n \frac{1}{2}\bZ/\bZ \rightarrow \frac{1}{2}\bZ/\bZ  \big)
\]
where $\divisor(M)=\frac{1}{2}\bZ/\bZ $ and $\Div(M)=0$.

\section{Preliminaries and reminder on global Galois cohomology}  \label{sec:prelim}


\subsection{Tate--Shafarevich and Selmer groups} \label{sec:tassel}

Let $k$ be an algebraic number field. For a discrete $\Gal_k$-module $M$ we set 
\[
\Sha^i(k,M) = \ker\Big(\rH^i(k,M) \xrightarrow{\res_v} \prod_v \rH^i(k_v,M)\Big)
\]
with the restriction maps $\res_v$ induced by the embedding $k \inj k_v$, and the product ranges over all places $v$ of $k$.  The Tate--Shafarevich group $\Sha(A/k)$ for an abelian variety $A$ over $k$ is defined\footnote{The traditional definition of $\Sha(A/k)$ is indeed equivalent due to the subtle equality
\[
\rH^1(k_v,A) = \rH^1(k_v,A(k_v^\alg)) = \rH^1(k_v,A(k^\alg)).
\]
} 
as $\Sha^1(k,A) = \Sha^1(k,A(k^\alg))$, and in particular it is a torsion group. This implies that  the Bashmakov problem can be considered one prime at a time.  We  will now concentrate on the $p$-primary part of the Bashmakov problem, i.e., analyzing the intersection of $\Sha(A/k)_{p^\infty}$ and the maximal divisible subgroup of the Weil-Ch\^atelet group.

Let $p$ be a prime number. The $p^n$-torsion Selmer group of $A$ is defined as
\[
\rH^1_\Sel(k,A_{p^n}) = \ker\Big(\rH^1(k,A_{p^n}) \to \prod_v \rH^1(k_v,A)\Big)
\]
with $v$ ranging over all places of $k$. A quick diagram chase with the cohomology sequence of the Kummer sequence 
\[
0 \to A_{p^n} \to A \xrightarrow{p^n\cdot} A \to 0
\]
over $k$ and all the localisations $k_v$ yields the fundamental short exact sequence
\begin{equation} \label{eq:selmerdef}
0 \to A(k)/p^nA(k) \xrightarrow{\delta_\kum} \rH^1_\Sel(k,A_{p^n}) \to \Sha(A/k)_{p^n} \to 0.
\end{equation}
It is known that $ \rH^1_\Sel(k,A_{p^n})$ is a finite group. 


\subsection{Generalized Selmer groups}
The Selmer group $\rH^1_\Sel(k,A_{p^n})$ is a global $\rH^1$ with local Selmer-conditions at every place. A \textbf{generalized Selmer group} $\rH^1_\bL(k,M)$ for a discrete finitely generated $\Gal_k$-module $M$ occurs by imposing \textbf{local conditions} $\bL_v \subset  \rH^1(k_v,M)$ as follows
\[
\rH^1_\bL(k,M) = \ker\Big(\rH^1(k,M) \to \prod_v \rH^1(k_v,M)/\bL_v\Big),
\]
such that the local conditions $\bL_v$ agree for almost all $v$ with the \textbf{unramified local cohomology}  
\[
\rH^1_{\nr}(k_v,M) = \infl\big(\rH^1(\kappa(v),M^{\rI_v})\big) = \ker\big(\rH^1(k_v,M) \to \rH^1(\rI_v,M)\big) .
\]
Here $\kappa(v)$ is the residue field at $v$ and $\rI_v \subset \Gal_{k_v}$ is the inertia subgroup. We do not bother to define $\rH^1_{\nr}$ for infinite places $v$ since there are only finitely many of them. For a textbook reference we refer to \cite{nsw} (8.7.8).

Observe that in the case when $M=A_{p^n}$, it is known that the image of the boundary map $\delta_{\kum}$ of the Kummer sequence 
\[
\Sel_v = \delta_\kum\big(A(k_v)/p^nA(k_v)\big) \subset \rH^1(k_v,A_{p^n})
\]
coincides with $\rH^1_{\nr}(k_v,A_{p^n})$ at places of good reduction with residue characteristic distinct from $p$ (see \cite{Milne} Theorem 2.6). Hence for $M=A_{p^n}$ with $\bL_v = \Sel_v$ we find back our original definition of the Selmer group.
 
\smallskip

Let $\bL=(\bL_v)$ be local conditions for the $\Gal_k$-module $M$, and let $Q$ be a finite set of places of $k$. Then we set $\bL_Q$ (resp.\ $\bL^Q$) for the local conditions \textbf{$\bL$ but free at $Q$} (resp.\ \textbf{$\bL$ but trivial at $Q$}) which agree with $\bL$ at all places $v \notin Q$ and with $\bL_{Q,v} = \rH^1(k_v,M)$ (resp.\ $\bL_v^Q = 0$) for all $v \in Q$.
We shall mainly be working with the following generalized Selmer groups. Let $Q$ be a finite set of finite places. The \textbf{Selmer group free at $Q$} is defined as 
\[
\rH^1_{\Sel_Q}(k,A_{p^n}) = \ker\Big(\rH^1(k,A_{p^n}) \to \prod_{v \not\in Q} \rH^1(k_v,A_{p^n})/\Sel_v\Big)
\]
and  the \textbf{Selmer group trivial at $Q$} is defined as 
\[
\rH^1_{\Sel^Q}(k,A_{p^n}) = \ker\Big(\rH^1(k,A_{p^n}) \to \prod_{v \not\in Q} \rH^1(k_v,A_{p^n})/\Sel_v \times \prod_{v \in Q} \rH^1(k_v,A_{p^n}) \Big).
\]
In the case when $Q$ is the set of primes of $k$ dividing $p$, we use $\rH^1_{\Sel^p}(k,A_{p^n})$ (resp. $\rH^1_{\Sel_p}(k,A_{p^n})$) to denote the Selmer groups $\rH^1_{\Sel^Q}(k,A_{p^n})$ (resp. $\rH^1_{\Sel_Q}(k,A_{p^n})$).
Obviously we have inclusions 
\[
\rH^1_{\Sel^Q}(k,A_{p^n}) \subseteq \rH^1_{\Sel}(k,A_{p^n}) \subseteq \rH^1_{\Sel_Q}(k,A_{p^n}).
\]


\subsection{Dual conditions and Euler characteristic formula}\label{sect:dualCond}
Let $\bL = (\bL_v)$ be a collection of local conditions for a discrete finitely generated $\Gal_k$-module $M$. The dual local conditions $\bL^\ast = (\bL_v^\ast)$ are defined as the orthogonal complements $\bL_v^\ast = \bL_v^\perp$ with respect to the local Tate-duality pairing. Hence $\bL^\ast$ is a collection of local conditions for the dual $\Gal_k$-module 
\[
M^D = \Hom(M,\bQ/\bZ(1)).
\]
We know that 
\[
\rH^1_{\nr}(k_v,M)=\rH^1_{\nr}(k_v,M^D)^\perp
\]
for almost all $v$, see \cite{Milne} Theorem~2.6, and
consequently,  if $\rH^1_\bL(k,M)$  is a generalized Selmer group then so is $\rH^1_{\bL^\ast}(k,M^D)$. 
If $M=A_{p^n}$ for an abelian variety $A/k$, so that $M^D = A^t_{p^n}$, where $A^t$ is the dual abelian variety, then the Selmer condition is self-dual: $\Sel^\ast = \Sel$. 

\smallskip

In order to relate sizes of these generalized Selmer groups, one adapts the Tate--Poitou exact sequence 
as in the proof of  \cite{nsw} Theorem (8.7.9).  We denote by $(-)^\vee$ the Pontrjagin dual $\Hom(-,\bQ/\bZ)$ and get the exact sequence
\begin{equation} \label{eq:tatepoitou}
0 \to \frac{\rH^2(k_S/k,M^D)^\vee}{\prod_{v \in S} \rH^0(k_v,M)} \to \rH^1_{\bL}(k,M) \to \prod_{v \in S} \bL_v \to \rH^1(k_S/k,M^D)^\vee \to \rH^1_{\bL^\ast}(k,M^D)^\vee \to 0
\end{equation}
where $S$ is a finite set of primes of $k$ which includes all the archimedean primes, all primes that ramify in the extension $k(M)/k$, all the primes for which $\bL_v \neq \rH^1_{\nr}(k_v,M)$, as well as all the primes dividing the order of $M$. Here, as usual, the field $k_S$ is the maximal extension of $k$ unramified outside $S$.

For two sets of local conditions $\bL_0 \subseteq \bL$, meaning $\bL_{0,v} \subseteq \bL_v$ for all $v$,  we use the natural map of \eqref{eq:tatepoitou} for $\bL_0$ to the sequence for $\bL$ with the same set $S$ adapted to both local conditions,
\[ 
\xymatrix@M+0.1ex@R-2ex@C-2.8ex{
0 \ar[r] & \frac{\rH^2(k_S/k,M^D)^\vee}{\prod_{v \in S} \rH^0(k_v,M)} \ar@{=}[d] \ar[r] &  \rH^1_{\bL_0}(k,M) \ar[r] \ar@{^(->}[d] &  \prod_{v \in S} \bL_{0,v} \ar@{^(->}[d] \ar[r] &  \rH^1(k_S/k,M^D)^\vee \ar@{=}[d] \ar[r] & \rH^1_{\bL_0^\ast}(k,M^D)^\vee \ar@{->>}[d] \ar[r] & 0
\\
0 \ar[r] & \frac{\rH^2(k_S/k,M^D)^\vee}{\prod_{v \in S} \rH^0(k_v,M)} \ar[r] &  \rH^1_{\bL}(k,M) \ar[r] &  \prod_{v \in S} \bL_v \ar[r] &  \rH^1(k_S/k,M^D)^\vee \ar[r] & \rH^1_{\bL^\ast}(k,M^D)^\vee \ar[r] & 0.
}
\] 
A diagram chase leads to
\[ 
\xymatrix@M+1ex@R-2ex@C-2.8ex{
0 \ar[r] & \frac{ \rH^1_{\bL_0}(k,M)}{\rH^2(k_S/k,M^D)^\vee}  \ar[r] \ar@{^(->}[d] &  \prod_{v \in S} \bL_{0,v} \ar@{^(->}[d] \ar[r] &  \ker\Big(\rH^1(k_S/k,M^D)^\vee \to \rH^1_{\bL_0^\ast}(k,M^D)^\vee \Big) \ar@{^(->}[d] \ar[r] & 0
\\
0 \ar[r] & \frac{ \rH^1_{\bL}(k,M)}{\rH^2(k_S/k,M^D)^\vee}  \ar[r] &  \prod_{v \in S} \bL_v  \ar[r] & \ker\Big( \rH^1(k_S/k,M^D)^\vee \to \rH^1_{\bL^\ast}(k,M^D)^\vee\Big) \ar[r] & 0.
}
\] 
Applying the snake lemma twice we find thus a short exact sequence
\begin{equation} \label{eq:selmerduality}
0 \to \frac{\rH^1_{\bL}(k,M)}{\rH^1_{\bL_0}(k,M)}  \to \prod_{v} \bL_v/\bL_{0,v} \to
 \left(\frac{\rH^1_{\bL_0^\ast}(k,M^D)}{\rH^1_{\bL^\ast}(k,M^D)}\right)^\vee \to 0
\end{equation}
where we can afford to take the product over all places $v$ since $\bL_{0,v} = \bL_v$ for $v \notin S$ by assumption.
We will apply \eqref{eq:selmerduality} in the following special case which will help estimating sizes of generalized Selmer groups later.

\begin{prop} \label{prop:epnumerics}
Let $M$ be a finite self-dual $\Gal_k$-module and let $\bL=\bL^\ast$ be a system of local conditions for $M$ that is self dual with respect to the identification $M \cong M^D$. Let $Q$ be a finite set of places. Then we have the numerical Euler-Poincar\'e characteristic formula
\[
\frac{| \rH^1_{\bL_Q}(k,M)|}{|\rH^1_{\bL^Q}(k,M)|} = \prod_{v \in Q}  |\bL_v|.
\]
\end{prop}
\begin{proof}
We count for $\bL \subseteq \bL_Q$ by $\rH^1_{\bL}(k,M) = \rH^1_{\bL^\ast}(k,M^D)$ and 
by \eqref{eq:selmerduality} 
\[
\frac{| \rH^1_{\bL_Q}(k,M)|}{|\rH^1_{\bL^Q}(k,M)|} = \frac{|\rH^1_{\bL_Q}(k,M)|}{|\rH^1_{\bL}(k,M)|} \cdot \frac{|\rH^1_{\bL^\ast}(k,M^D)|}{|\rH^1_{\bL_Q^\ast}(k,M^D)|} = \prod_{v \in Q}  |\rH^1(k_v,M)/\bL_v| = \prod_{v \in Q}  |\bL_v|,
\]
where the last equality is due to 
\[
|\rH^1(k_v,M)/\bL_v| = |\Hom(\bL_v^\ast,\bQ/\bZ)| = |\bL_v^\ast| = |\bL_v|
\]
because of the self duality of $\bL$.
\end{proof}


\section{Translation between Galois cohomology and \'etale cohomology} \label{sec:ec}


\subsection{Compactly supported \'etale cohomology}

In this section we shall not be concerned with the effect of real places by assuming there are none, or that we deal with $p \not=2$. Otherwise, there are expositions of the necessary modifications after Artin--Verdier by introducing the Woods-Hole site,
see for example \cite{zink:appendix}. 

Set $X=\Spec(\fo_k)$ 
for the ring of integers $\fo_k$ of an algebraic number field $k$. Let $j: U \subset X$ be a dense Zariski open. Then compactly supported cohomology of a constructible torsion sheaf $\cF$ on the small \'etale site $U_\et$ is defined as 
\[
\rH^i_{\rc}(U,\cF) = \rH^i(X,j_! \cF).
\]
Restriction induces the  natural "forget support" map $\rH^i_{\rc}(U,\cF) \to \rH^i(U,\cF)$.


\subsection{\'Etale cohomology and the relation to generalized Selmer groups}

Let $\cM$ be a locally constant constructible sheaf on $U_\et$ with generic fiber the finite  $\Gal_k$-module $M$. Part of the  localisation sequence for an open $V \subset U$ reads 
\begin{equation} \label{eq:localisation}
\bigoplus_{v \in U \setminus V} \rH^1_v(U,\cM) \to \rH^1(U,\cM) \to \rH^1(V,\cM) \to \bigoplus_{v \in U \setminus V}  \rH^2_v(U,\cM).
\end{equation}
As $\rH^1_v(U,\cM)=0$ for locally constant sheaves $\cM$, we find that 
\[
\rH^1(U,\cM) \inj \rH^1(k,M) = \varinjlim_V \rH^1(V, \cM)
\]
is injective. Let $\bL = (\bL_v)$ be a collection of local conditions for $M$. We define the \'etale cohomology with local conditions as 
\begin{equation} \label{eq:defigeneralizedselmer}
\rH^1_\bL(U,\cM) = \ker\Big(\rH^1(U,\cM) \to \prod_v \rH^1(k_v,M)/\bL_v\Big),
\end{equation}
which agrees with  $\rH^1(U,\cM) \cap \rH^1_\bL(k,M)$ inside $\rH^1(k,M)$.
In particular, we have \'etale cohomology $\rH^1_{\Sel}(U,A_{p^n})$ with Selmer condition for a Zariski open $U \subset X$ with good reduction of $A$ over  $U$ and a prime number $p$ invertible on $U$.

\begin{prop}  \label{prop:comparegeneralizedselmer}
 If $\bL_v = \rH^1_{\nr}(k_v,M)$ for all $v \in U$ then 
 \begin{equation} \label{eq:comparegeneralizedselmer}
 \rH^1_\bL(U,\cM) = \rH^1_\bL(k,M).
 \end{equation}
If in addition $\bL_v = \rH^1(k_v,M)$  for all $v \not\in U$, then 
\[
\rH^1(U,\cM) = \rH^1_\bL(k,M),
\]
whereas if in addition $\bL_v = 0$  for all $v \not\in U$, then 
\[
\im\big(\rH^1_{\rc}(U,\cM) \to \rH^1(U,\cM)\big) = \rH^1_\bL(k,M).
\]
\end{prop}

\begin{proof}
By definition of $\rH^1_\bL(U,\cM)$ it suffices to show that $\rH^1_\bL(k,M) \subset \rH^1(U,\cM)$.
Any class $\alpha \in \rH^1(k,M)$ lies in $\rH^1(V,\cM)$ for small enough open $V \subset U$. The claim follows from (\ref{eq:localisation}) because the image of an $\alpha \in \rH^1_\bL(k,M)$ vanishes in  
\[
\rH^2_v(U,\cM) = \rH^1(k_v,M)/\rH^1_{\nr}(k_v,M).
\]

The additional claims follow from (\ref{eq:comparegeneralizedselmer}), the definition (\ref{eq:defigeneralizedselmer}), and the exact sequence
\begin{equation} \label{eq:compactlocalisation}
 \bigoplus_{v \in X \setminus U} \rH^{i-1}(k_v,M) \to \rH^i_{\rc}(U,\cM) \to \rH^i(U,\cM) \xrightarrow{\res_v}  \bigoplus_{v \in X \setminus U} \rH^i(k_v,M)
\end{equation}
for $i \geq 0$. These exact sequences arise from the localisation sequence for $U\subset X$ and the sheaf $j_!\cM$ because by excision
$\rH^{i+1}_v(X,j_!\cM)=\rH^i(k_v,M)$.
\end{proof}

\begin{cor} \label{cor:comparegeneralizedselmer}
Let $A/k$ be an abelian variety with  good reduction over $U$ and let $p$ be a rational prime invertible on $U$. With $Q$ equal to the set of finite places $v \not\in U$ we have 
\begin{eqnarray}
\rH^1_{\Sel_Q}(k,A_{p^n}) & = &  \rH^1(U,A_{p^n}) \\
\rH^1_{\Sel^Q}(k,A_{p^n}) & = & \im\big(\rH^1_{\rc}(U,A_{p^n}) \to \rH^1(U,A_{p^n})\big) 
\end{eqnarray}
\end{cor}
\begin{proof}
By assumption for $v \in U$ we have $\Sel_v = \rH^1_{\nr}(k_v,A_{p^n})$, so the corollary follows from Proposition~\ref{prop:comparegeneralizedselmer}.
\end{proof}

\begin{cor} \label{cor:ShaInGlobal}
Let $U \subset X$ be a Zariski open where the abelian variety $A/k$ has good reduction and $p$ is invertible. By abuse of notation we also denote by $A \to U$ the smooth model of $A/k$. Then for an open $V \subset U$ the restriction map $\rH^1(U,A) \to \rH^1(V,A)$ is injective and 
\[
 \Sha(A/k)_{p^\infty} \subseteq \rH^1_{\divisor}(k,A)_{p^\infty} \subseteq  \rH^1(U, A) \subseteq \rH^1(k,A) = \varinjlim_{V \subset U} \rH^1(V, A).
\]
\end{cor}
\begin{proof}
An element of $\rH^1(U,A)$ is represented by a principal homogeneous space $W$ under $A$  over $U$, which is trivial if and only $W(U)$ is nonempty. The valuative criterion of properness yields $W(U) = W(V)$ which proves the claim on injectivity. 

Now $\rH^1_{\divisor}(k,A)_{p^n}$ consists of classes that are trivial outside all $v \mid p$, see 
\eqref{eq:pparth1localdiv}, hence is contained in the image of 
\[
\rH^1_{\Sel_p}(U,A_{p^n}) = \rH^1_{\Sel_p}(k,A_{p^n}) \to \rH^1(k,A)
\]
which maps to $\rH^1(U,A)$.
\end{proof}


\section{Review of Bashmakov's results: when the pro-\texorpdfstring{$p$}{p} fundamental group is free}  \label{sec:bashmakov}

In \cite{bashmakov:div} Bashmakov proves that 
\[
\Sha(A/k)_{p^\infty} \subseteq \Div(\rH^1(k,A))
\]
for abelian varieties $A/k$ with very special properties. First, the number field $k$ is required to satisfy  the following conditions:
\begin{enumerate}
\item[(i)] $k$ contains the $p$th roots of unity $\zeta_p$, 
\item[(ii)] there is a unique $\fp \mid (p)$  in $k/\bQ$,
\item[(iii)] and the class number of $k$ is prime to $p$. 
\end{enumerate}
It is a well known result of Shafarevich \cite{shafarevich},
 that conditions (i)--(iii) imply that the complement 
\[
V = \Spec(\fo_k) \setminus \{\fp\}
\]
has a free pro-$p$ group as its maximal pro-$p$ quotient $\pi_1(V)^{\text{pro}-p}$ of its fundamental group, see \cite{nsw} Cor 10.7.14. Moreover, the abelian variety $A/k$ has to satisfy:
\begin{enumerate}
\item[(iv)] $A/k$ has good reduction above $V$,
\item[(v)] the action of $\Gal_k$ on $A_p$ factors over a $p$-group.
\end{enumerate}
In \cite{bashmakov:div}, Bashmakov does not require condition (v), although his argument at the very end does depend on it. Namely, under these assumptions  we compute by \cite{zink:appendix} Proposition 3.3.1
\[
\rH^2(V,A_{p^n}) = \rH^2(\pi_1(V),A_{p^n}).
\]
Due to (v) and a theorem of Neumann \cite{neumann:pclosed}, see \cite{nsw} Cor 10.4.8, we find that
\begin{equation} \label{eq:withconditionv}
\rH^2(\pi_1(V),A_{p^n}) = \rH^2(\pi_1(V)^{\text{pro}-p},A_{p^n}) = 0 
\end{equation}
which vanishes in view of the freeness of $\pi_1(V)^{\text{pro}-p}$ recalled above.

\begin{prop} \label{prop:bashmakovfreepi1}
Let $p$ be a regular prime number. Let $k$ be a subfield of $\bQ(\zeta_p)$ and let $A/k$ be an abelian variety with good reduction away from the prime above $p$. If $\Gal_{\bQ(\zeta_p)}$ acts via a $p$-group on $A_p$, then $\rH^1_{\divisor}(k,A)_{p^\infty}$ agrees with the $p$-part of $\Div(\rH^1(k,A))$.
\end{prop}

\begin{proof} 
We set $U=\Spec(\fo_k[1/p])$ and write by abuse of notation $A/U$ for the smooth model of $A/k$ over $U$.  
By Corollary~\ref{cor:ShaInGlobal} we have the second and third inclusion in 
\[
\Div(\rH^1(k,A))_{p^\infty} \subseteq \rH^1_{\divisor}(k,A)_{p^\infty} \subseteq \rH^1(U,A)_{p^\infty}  \subseteq \rH^1(k,A)_{p^\infty}.
\]
Thus, and since $\rH^1(U,A)_{p^\infty}$ is a finitely cogenerated abelian torsion group, we find
\[
\Div(\rH^1(k,A))_{p^\infty} = \Div(\rH^1(U,A)_{p^\infty}) = \divisor(\rH^1(U,A)_{p^\infty}).
\]
It therefore suffices to show that every element of $\rH^1(U,A)_{p^\infty}$ is $p$-divisible in $\rH^1(U,A)_{p^\infty}$. The Kummer  exact sequence 
\[
0 \to A_{p^r} \to A \xrightarrow{p^r \cdot} A \to 0
\]
on $U_\et$ yields a short exact sequence
\[
\rH^1(U, A) \xrightarrow{p^r \cdot}  \rH^1(U, A) \xrightarrow{\delta_r}  \rH^2(U, A_{p^r}) 
\]
and the task left is showing $\delta_r=0$.
Since $\bQ(\zeta_p)/k$ is of order prime to $p$, restriction via the finite \'etale 
\[
V=\Spec(\bZ[\zeta_p,1/p]) \to U
\]
is injective by the corestriction argument and we have by \eqref{eq:withconditionv}
\[
\rH^2(U,A_{p^n}) \inj \rH^2(V,A_{p^n}) =  \rH^2(\pi_1(V)^{\text{pro}-p},A_{p^n}) = 0,
\]
completing the proof.
\end{proof}

\begin{ex}
Examples of abelian varieties above $k=\bQ(\zeta_p)$ with respect to a regular prime $p$ and which satisfy the constraints imposed by Bashmakov (including condition (v)) are given by the Jacobian of the Fermat curve $X^p + Y^p = 1$, see for example \cite{andersonihara:proell}.  
\end{ex}

However, the natural second family of examples to consider, namely the Jacobian of the modular curve $X_0(p)$ for $p=11$ or $p \geq 17$, fail condition (v). Indeed\footnote{We thank K.~Ribet for information on the Galois representation of the modular Jacobian.}, no prime $\fm$ of the Hecke algebra $\bT$ above $p$ is Eisenstein (those divide $p-1$) and the corresponding representation $\rho_\fm : \Gal_{\bQ} \to \GL_2(\bT/\fm)$ is irreducible by Mazur \cite{mazur:eisenstein}. Thus $\rho_\fm$ is surjective due to  \cite{ribet:sstgalois}.


\section{Finite subgroups of \texorpdfstring{$\GL_2$}{GL-zwei}}  \label{sec:group}
The purpose of this section is to provide the necessary statements from finite group theory.

\begin{lem} \label{lem:center} Let $W$ be a finite dimensional $\bF_p$-vector space, and let $G \subset \GL(W)$ be a subgroup which intersects the center $\bF_p^\ast \cong Z \subset \GL(W)$ nontrivially. Then the following holds.
\begin{enumerate}
\item $W$ and the adjoint representation $\End(W)$ have no common irreducible factor.
\item $\rH^1(G,W) = 0$.
\end{enumerate}
\end{lem}
\begin{proof}
(1) The group $H=G \cap Z$ of homotheties in $G$ acts trivially on every irreducible factor of $\End(W)$ and faithfully on every irreducible factor of $W$. Hence none of them can occur in both $W$ and $\End(W)$.

(2) The inflation/restriction sequence for $H \unlhd G$ reads 
\[
0 \to \rH^1(G/H,W^H) \to \rH^1(G,W) \to \rH^1(H,W)^{G/H}.
\]
Since $H$ was assumed nontrivial and is necessarily of order prime to $p$, both $W^H$ and $ \rH^1(H,W)$ vanish, and consequently also $ \rH^1(G,W) = 0$. 
\end{proof}

\smallskip

In the $2$-dimensional case we have a further criterion.

\begin{thm} \label{thm:groupcrit}
Let $V$ be a vector space over $\bF_p$ of dimension $2$ with $p \geq 3$, and let $G \subseteq \GL(V)$ be a subgroup such that 
the $G$-module $V$ is irreducible.
Then the following holds:
\begin{enumerate}
\item[(1)] $\rH^1(G,V) = 0$,
\end{enumerate}
and if, moroeover, $G$ is not conjugate to a subgroup of $S_3 \subseteq \GL(V)$ with respect to a $2$-dimensional irreducible representation of $S_3$, then:
\begin{enumerate}
\item[(2)] $V$ and $\End(V)$ have no common irreducible factor as $G$-modules.
\end{enumerate}
\end{thm}

\begin{proof}
We first recall the 
classification of subgroups of $\GL_2(\bF_p)$, see \S2 of \cite{Serre}.
Let $\ov{G}$ be the image of $G$ under the natural map $\GL_2(\bF_p) \to \PGL_2(\bF_p)$. Then one of the following  holds.
\begin{enumerate}
\item[(a)] $p \mid \#G$ and $G$ is contained in a Borel $B \subset \GL_2(\bF_p)$.
\item[(b)] $p \mid \#G$ and $G$ contains $\SL_2(\bF_p)$.
\item[(c)] $p \nmid \#G$ and $G$ is contained in a normalizer of a split torus.
\item[(d)] $p \nmid \#G$ and $G$ is contained in a normalizer of a non-split torus.
\item[(e)] $p \nmid \#G$ and $\ov{G}$ is isomorphic to $A_4$, $S_4$, or $A_5$.
\end{enumerate}

\textit{The case $p \mid \#G$.} 
As by assumption $G$ is not contained in a Borel, we conclude by the above list  that $G$ contains $\SL(V)$. Now since $p\geq 3$,  the group $G$ necessarily meets the center of $\GL(V)$ nontrivially, so that we conclude by Lemma~\ref{lem:center}. 

\textit{The case $p \nmid \#G$.} 
In this case $\rH^1(G,V) = 0$ holds trivially, so we are done with assertion (1).  If $G$ belongs to case (e) of the above list, then Lemma~\ref{lem:exceptionalcase} below shows that the intersection with the center of $\GL(V)$ is nontrivial and we conclude by Lemma~\ref{lem:center}.

It remains to discuss the case that $p \nmid \#G$ and $G$ is contained in the normalizer $N = C \rtimes \bZ/2\bZ$ of a torus $C \subset \GL(V)$ such that $V$ is an irreducible $G$-module. In order to understand the $G$-action on $\End(V)$ we identify $V  =
\bF_p[\alpha]$ with a separable quadratic $\bF_p$-algebra with $\Aut(\bF_p[\alpha]/\bF_p) = \bZ/2\bZ$ generated by 
$F \in \End(V)$ such that 
\[
N = \bF_p[\alpha]^\ast \rtimes \langle F \rangle \subset \GL(V).
\]
We find an isomorphism
\[
\bF_p[\alpha] \oplus \bF_p[\alpha] \cdot F \xrightarrow{\sim} \End(V) 
\]
\[
a + b\cdot F \mapsto x \mapsto ax + bF(x)
\]
as $N$-modules, where $N$ acts on the first summand $V_1 = \bF_p[\alpha]$ through the quotient 
\[
 \bF_p[\alpha]^\ast \rtimes \langle F \rangle \surj \langle F \rangle = \Aut(\bF_p[\alpha]/\bF_p) \subset \GL(V),
\]
and, by a straight forward calculation left to the reader, on $V_2=\bF_p[\alpha] \cdot F$ through the endomorphism 
\[
 \bF_p[\alpha]^\ast \rtimes \langle F \rangle \to    \bF_p[\alpha]^\ast \rtimes \langle F \rangle 
\]
which maps  $F$ to $F$ and $\lambda \in \bF_p[\alpha]^\ast$ to $\lambda/F(\lambda)$, followed by the tautological action on 
$V = \bF_p[\alpha]$ identified with $\bF_p[\alpha] \cdot F$ by formally multiplying with $F$. 
As $p \nmid \#G$ the representation theory of $G$  in $\bF_p$-vector spaces is semissimple and assertion (2) can only fail if $V$ is isomorphic to $V_1$ or $V_2$.

Let $H = G \cap C$ be the intersection with the torus, hence a subgroup in $G$ of index $\leq 2$. Because $V$ is not a reducible $G$-module and $p \geq 3$ we have $\# G \geq 3$ and therefore $H \not=1$.  As $H$ acts trivial on $V_1$  
we have $V \not\cong V_1$ as $G$-modules. It remains to exclude $V \cong V_2$ as $G$-modules.

The representation $V \otimes \bF_{p^2}$ regarded as an $H$-module decomposes as a
sum of characters $\chi_1 \oplus \chi_2$ of $H$. The representation $V \cdot F \subset \End(V)$ decomposes after scalar extension to $\bF_{p^2}$ as $H$-module as $\chi_1 \chi_2^{-1} \oplus \chi_2 \chi_1^{-1}$. Comparing the two, we find either $\chi_1 = \one = \chi_2$, whence $H=1$ contradicting the irreduciblity of $V$ as a $G$-module. Or, $\chi_1  =\chi_2^2 \not=1$ and $\chi_2 = \chi_1^2\not=1$ which means $\chi_1$ and $\chi_2$ are of order $3$ and determine each other. In this case $H \cong \bZ/3\bZ$  and acts on $V$ noncentrally and without a common fixed vector. In any case, split or non-split, the subgroup $H \subset C$ is normal but not central in $N$. Hence  either $H = G$ and $G$ can be embedded in a subgroup $S_3 \subseteq \GL(V)$, or $H \unlhd G$ of index $2$ and $G \cong S_3$ itself. In any case, these subgroups $G$ are excluded by assumption.
This completes the proof of Theorem~\ref{thm:groupcrit}.
\end{proof}

\begin{lem} \label{lem:exceptionalcase}
Let $G$ be a subgroup of $\GL_2(\bF_p)$ such that  $p \nmid \#G$ and in the notation above $\ov{G}$ is isomorphic to $A_4$, $S_4$, or $A_5$, i.e., in case (e) of the above list. Then $G$ meets the center $Z$ of $\GL_2(\bF_p)$ nontrivially.
\end{lem}
\begin{proof}
It suffices to discuss $\ov{G} = A_4$. If $G \cap Z = 1$, then we have a copy $A_4 \subseteq \GL_2(\bF_p)$ and $p >2$. The $2$-Sylow subgroup $V_4 \cong \bZ/2\bZ \times \bZ/2\bZ$ of $A_4$ has a completely reducible representation theory already rationally over $\bF_p$ as we may produce enough projectors already  over $\bF_p$. Hence $V_4$ is contained in a split torus $C = \bF_p^\ast \times \bF_p^\ast$ and must agree with the $2$-torsion of $C$.   Thus $V_4$ already contains the central  element $-1$, a contradiction.
\end{proof}


\section{By induction to  the structure of generalized Selmer groups} \label{sec:induction}


\subsection{The Selmer splitting field}   \label{sec:Selmersplitting}
Let $A/k$ be an abelian variety. 
We set 
\[
\rH^1_{\Sel}(k(A_p)/k,A_p)
\]
for the intersection of $\rH^1(k(A_p)/k,A_p)$ under inflation with $\rH^1_{\Sel}(k, A_p)$ in $\rH^1(k,A_p)$.
Then we have the following commutative diagram with exact rows 
\[ 
\xymatrix@M+1ex@R-2ex{
0 \ar[r] &  \rH^1_{\Sel}(k(A_p)/k,A_p) \ar@{}[d]|{{\rotatebox{-90}{$\subseteq$}}} \ar[r]^(0.55){\infl} & \rH^1_{\Sel}(k, A_p) \ar[r]^(0.3){\res}  \ar@{}[d]|{{\rotatebox{-90}{$\subseteq$}}} & \Hom_{\Gal(k(A_p)/k)}(\Gal_{k(A_p)}^\ab \otimes \bF_p,A_p)  \ar@{=}[d] \\
0 \ar[r] &  \rH^1(k(A_p)/k,A_p) \ar[r]^(0.55){\infl} & \rH^1(k, A_p) \ar[r]^(0.3){\res} & \rH^1(k(A_p),A_p)^{\Gal(k(A_p)/k)} 
} 
\] 
The exactness of the top row follows from the exactness of the bottom row. The restriction map defines a canonical continuous $\Gal_k$-equivariant pairing 
\[
 \rH^1_{\Sel}(k, A_p) \times \big(\Gal_{k(A_p)}^\ab \otimes \bF_p\big) \to A_p.
\]
Let $\Gal_{k(A_p)}^\ab \otimes \bF_p \surj M$ be the continuous finite  quotient by the right kernel of the pairing. Then the restriction  map factors as 
\[
\rH^1_{\Sel}(k, A_p) \to  \Hom_{\Gal(k(A_p)/k)}(M,A_p)  \subset   \Hom_{\Gal(k(A_p)/k)}(\Gal_{k(A_p)}^\ab \otimes \bF_p,A_p) .
\]
The quotient $M$ corresponds to a finite Galois extension $L/k(A_p)$,  the \textbf{Selmer splitting field} of $A$ with respect to the prime $p$, more precisely $M=\Gal(L/k(A_p))$. Since $M$ is a quotient as $\Gal(k(A_p)/k)$-module, the field $L$ is in fact Galois over $k$ and $\Gal(k(A_p)/k)$ acts on $M$ by conjugation after lifting under the quotient map $\Gal(L/k) \surj \Gal(k(A_p)/k)$. 

\begin{lem} \label{lem:selmerslpittingfield}
Let $A/k$ be an abelian variety, and let $L$ be the Selmer splitting field with respect to $p$. Then  the following holds.
\begin{enumerate}
\item The following sequence is exact:
\begin{equation} \label{eq:selmersplittingfield}
0 \to \rH^1_{\Sel}(k(A_p)/k,A_p) \to \rH^1_{\Sel}(k, A_p) \to  \Hom_{\Gal(k(A_p)/k)}(\Gal(L/k(A_p)),A_p)
\end{equation}
\item Every irreducible $\Gal(k(A_p)/k)$-module subquotient of $\Gal(L/k(A_p))$ is isomorphic to an irreducible subquotient of $A_p$.
\end{enumerate}
\end{lem}
\begin{proof}
(1) is obvious by the definition of $L$. For (2) we note that the defining pairing yields an injective map
\begin{equation} \label{eq:structureM}
\Gal(L/k(A_p)) = M \inj \Hom( \rH^1_{\Sel}(k, A_p),A_p) \cong A_p \oplus \ldots \oplus A_p
\end{equation}
of  $\Gal(k(A_p)/k)$-modules, where $\rH^1_{\Sel}(k, A_p)$ carries trivial action.
\end{proof}

\subsection{The structure of some generalized Selmer groups} The following theorem will be ultimately applied to elliptic curves that are automatically principally polarized.

\begin{thm}  \label{thm:structureofgenselmer}
Let $A/k$ be a principally polarized abelian variety, and let $p$ be a prime number, such that
\begin{enumerate}
\item[(i)] $A_p(k) = 0$, 
\item[(ii)] $\rH^1(k(A_p)/k,A_p) = 0$.
\end{enumerate}
Let $Q$ be a finite set of finite primes of $k$ not dividing $p$, and fix $n \in \bN$ such that 
\begin{enumerate}
\item[(iii)] $A$ has good reduction at $v$ for all $v \in Q$;
\item[(iv)]  the set of Frobenius elements $\Frob_{w} \in \Gal(L/k(A_p))$ where $L$ is the Selmer splitting field of $A/k$ with respect to $p$ and $w$ denotes a prime of $k(A_p)$ dividing $v$, when $v$ ranges over $Q$,  generates $\Gal(L/k(A_p))$ as a $\Gal(k(A_p)/k)$-module;
\item[(v)] $A_{p^n}(k_v)$ is a free $\bZ/p^n\bZ$-module for all $v \in Q$.
\end{enumerate}
Then for all $m \leq n$ we have that
\begin{enumerate}
\item $\rH^1_{\Sel^Q}(k,A_{p^m}) = 0$,
\item $\rH^1_{\Sel_Q}(k,A_{p^m}) \cong \prod_{v \in Q} A_{p^m}(k_v)$ under a non-canonical group isomorphism.
\end{enumerate}
\end{thm}

\begin{proof}
\textit{Step 1:}  We first treat (1) for $m=1$.  We set $k_1 = k(A_p)$, and $k_{1,w}$ for the completion of $k(A_p)$ in $w$. Localization at $v$ respectively $w$ yields a commutative diagram
\[
\xymatrix@M+1ex@R-1ex@C+2ex{ \rH^1_{\Sel}(k, A_p) \ar[r]^(0.35){\res_{k_1/k}} \ar[d] & \Hom_{\Gal(k_1/k)}(\Gal(L/k_1),A_p) \ar[d]^{\ev_w} \\
\rH^1_{\nr}(k_v,A_p) \ar[r]^(0.45){\res_{k_{1,w}/k_v}} & \rH^1_{\nr}(k_{1,w},A_p) = A_p}
\]
with the evaluation map $\ev_w$ mapping a morphism $\ph : \Gal(L/k_1)  \to A_p$ to its value $\ph(\Frob_{w})$ at the Frobenius element of $w$. Assumption (iv), the sequence \eqref{eq:selmersplittingfield}, and assumption (ii) imply
\[
\rH^1_{\Sel^Q}(k,A_p) \subseteq \rH^1_{\Sel}(k(A_p)/k,A_p) \subseteq \rH^1(k(A_p)/k,A_p) = 0.
\]

\smallskip

\textit{Step 2:} We now show (1) by induction on $m$ terminating in $n$. As an abbreviation we set 
\[
\bL_{m,v} = \Sel^Q_v \subseteq \rH^1(k_v,A_{p^m})
\]
for the Selmer condition trivial at $Q$ for $A_{p^m}$-coefficients. Then 
 the following diagram is commutative and the rows are exact
\begin{equation} \label{eq:compatiblelocalconditions}
\xymatrix@M+1ex@R-2ex{   
 & {\bL_{m-1,v}} \ar[r] \ar[d] &  {\bL_{m,v}} \ar[r] \ar[d]  &  {\bL_{1,v}} \ar[d] \ar[r] & 0  \\
0 \ar[r] & \rH^1(k_v,A_{p^{m-1}})/\delta_\kum\big(A_p(k_v)\big) \ar[r] &  \rH^1(k_v,A_{p^m}) \ar[r]^{p^{m-1}\cdot} &  \rH^1(k_v,A_{p})
}
\end{equation}
The snake lemma applied to \eqref{eq:compatiblelocalconditions} yields that in the commutative diagram
\begin{equation} \label{eq:almostSel}
\xymatrix@M+1ex@R-2ex@C-2ex{ A_p(k) \ar[d] \ar[r]^(0.4){\delta_\kum} &  \rH^1(k,A_{p^{m-1}})  \ar[d] \ar[r] & \rH^1(k,A_{p^m}) \ar[r]  \ar[d] & \rH^1(k,A_p)  \ar[d] \\
   \prod_{v} A_p(k_v) \ar[r]^(0.4){\delta_\kum} &  \prod_{v } \bruch{\rH^1(k_v,A_{p^{m-1}})}{\bL_{m-1,v}} \ar[r] & \prod_{v} \bruch{\rH^1(k_v,A_{p^m})}{\bL_{m,v}} \ar[r] & \prod_{v} \bruch{\rH^1(k_v,A_p)}{\bL_{1,v}}
}
\end{equation} 
the bottom row is exact.  The map $\delta_\kum$ in the bottom row is the zero map: by assumption (v) when $v \in Q$, then $A_{p^m}(k_v)/A_{p^{m-1}}(k_v) \to A_p(k_v)$ is surjective for $m \leq n$, or, in general for $v \notin Q$, by comparing the boundary maps for the diagram 
\[
\xymatrix@M+1ex@R-2ex{
0 \ar[r] & A_{p^{m-1}} \ar@{=}[d] \ar[r] & A_{p^m} \ar@{^(->}[d]  \ar[r]^{p^{m-1} \cdot } & A_p \ar@{^(->}[d] \ar[r] & 0 \\
0 \ar[r] & A_{p^{m-1}} \ar[r] & A \ar[r]^{p^{m-1} \cdot } & A \ar[r] & 0 
}
\]
(It is the limitation of assumption (v) that forces the induction terminate at $n$). Again the snake lemma applied to 
\eqref{eq:almostSel} yields exactness of 
\[
\rH^1_{\Sel^Q}(k,A_{p^{m-1}})  \to \rH^1_{\Sel^Q}(k,A_{p^m}) \to \rH^1_{\Sel^Q}(k,A_{p}) 
\]
so that with Step 1 we deduce (1) by induction on $m$.

\smallskip

\textit{Step 3:} By (i) we have 
\[
\rH^1(k,A_{p^m}) =  \rH^1(k,A_{p^n})_{p^m},
\]
and it follows then from the definition that 
\[
\rH^1_{\Sel_Q}(k,A_{p^m}) =  \rH^1_{\Sel_Q}(k,A_{p^n})_{p^m}.
\]
is also  the exact $p^m$-torsion.

\smallskip

\textit{Step 4:} 
Since $A$ is principally polarized, the coefficients $A_{p^m}$ are self dual and the Selmer condition is self dual with respect to the identification $A_{p^m} = A_{p^m}^D$.

With the  Euler-Poincar\'e characteristic formula of Proposition~\ref{prop:epnumerics} and (1) we compute
\[
|\rH^1_{\Sel_Q}(k,A_{p^m})| = \frac{|\rH^1_{\Sel_Q}(k,A_{p^m})|}{|\rH^1_{\Sel^Q}(k,A_{p^m})|} = \prod_{v \in Q} |A(k_v)/p^mA(k_v)| = \prod_{v \in Q} |A_{p^m}(k_v)|.
\]
The last equation  $ |A(k_v)/p^mA(k_v)| = |A_{p^m}(k_v)|$ follows because by Mattuck--Tate the $p$-primary part of $A(k_v)$ is finite  since $v \nmid p$ for all $v \in Q$.

\smallskip

\textit{Step 5:} In order to show (2) it suffices to treat the case of $\rH^1_{\Sel_Q}(k,A_{p^n})$. By Steps 3 and 4 we conclude that $\rH^1_{\Sel_Q}(k,A_{p^n})$ has the same number of $p^m$-torsion elements for every $m \leq n$ as 
\[
\prod_{v \in Q} A_{p^n}(k_v),
\]
so that both groups are noncanonically isomorphic. This proves (2).
\end{proof}

\subsection{The case of elliptic curves} Based on the group theory of $\GL_2$ in Section~\ref{sec:group} we can show existence for auxiliary sets of primes $Q$ in Theorem~\ref{thm:structureofgenselmer} in the special case of elliptic curves. 

\begin{prop} \label{prop:structuregenselmerellipticcurve}
Let $E/k$ be an elliptic curve and let $p$ be an odd prime number such that $E_p$ is an irreducible $\Gal_k$-module and $\Gal(k(E_p)/k) \subseteq \GL(E_p)$ is not contained in a conjugate of $S_3 \subseteq \GL(E_p)$. Fix an $n\in \bN$. Then 
\begin{enumerate}
\item[(1)] $E_p(k) = 0$, 
\item[(2)] and $\rH^1(k(E_p)/k,E_p) = 0$.
\end{enumerate}
Moreover, we can find a finite set of primes $Q$ not dividing $p$ such that 
\begin{enumerate}
\item[(3)] $E$ has good reduction at $v$ for all $v \in Q$;
\item[(4)] the primes $v \in Q$ are completely split in $k(E_{p^n})/k$;
\item[(5)] the set of $\Frob_w \in \Gal(L/k(E_p))$ where $L$ is the Selmer splitting field of $E/k$ with respect to $p$ and $w$ denotes a prime of $k(E_p)$ dividing $v$  as $v$ varies though $Q$, generates $\Gal(L/k(E_p))$ as a  $\Gal(k(E_p)/k)$-module.
\end{enumerate}
In particular, for all $m \leq n$ we have
\begin{enumerate}
\item[(6)] $\rH^1_{\Sel_Q}(k,E_{p^m}) \cong (\bZ/p^m\bZ)^{2\cdot \#Q}$.
\end{enumerate}
\end{prop}

\begin{proof}
The subgroup $E_p(k) \subseteq E_p(k^\alg)$ is a $\Gal_k$-submodule and thus in view of the irreducibility assumption either all or nothing. In case of trivial $\Gal_k$-action, the module $E_p$ is not irreducible, so that we conclude (1).  Assertion (2) follows from Theorem~\ref{thm:groupcrit} (1).

\smallskip

We will now construct the set $Q$ of auxiliary primes. First we prove that $L$ and $k(E_{p^n})$ are linearly disjoint over $k(E_p)$. Indeed, let $K = L \cap k(E_{p^n})$ be their intersection and set $\ov{M} = \Gal(K/k(E_p))$ for the abelian Galois group over $k(E_p)$. Then, since $K/k$ is  Galois, the projection
\[
\Gal(L/k(E_p)) \surj \ov{M}
\]
is a surjection of $\Gal(k(E_p)/k)$-modules. It follows from Lemma~\ref{lem:selmerslpittingfield} that $\ov{M}$ has a composition series as $\Gal(k(E_p)/k)$-module consisting of irreducible subquotients of $E_p$.
On the other hand, the group $\Gal(k(E_{p^n})/k(E_p))$ is a subgroup of 
\[
N = \ker\big(\GL(E_{p^n}) \to \GL(E_p)\big).
\]
The group $N$ is solvable with abelian subquotients 
\[
N_m = \ker\big(\GL(E_{p^m}) \to \GL(E_{p^{m-1}})\big)
\]
that are canonically $\Gal(k(E_p)/k)$-modules by conjugation with lifts and as such isomorphic to the adjoint representation of $\Gal(k(E_p)/k)$ on $\End(E_p)$. 
Since, by Theorem~\ref{thm:groupcrit} (2),
$E_p$ and $\End(E_p)$ have no common irreducible $\Gal(k(E_p)/k)$-subquotient, we deduce that $\ov{M} = 0$ and $K = k(E_p)$, which means that $L$ and $k(E_{p^n})$ are linearly disjoint over $k(E_p)$.

The Chebotarev density theorem enables us to choose a finite set $Q$ of finite places $v \nmid p$ in the locus of good reduction of $E/k$, so that the Frobenius elements $\Frob_v$ for $v \in Q$ satisfy
\begin{enumerate}
\item[(i)] the image of $\Frob_v$ in $\Gal(k(E_{p^n})/k)$ is trivial,
\item[(ii)] the  images of $\Frob_v$ for $v \in Q$ generate $\Gal(L/k(E_p))$.
\end{enumerate}
The linear disjointness of $L$ and $k(E_{p^n})$ over $k(E_p)$ implies that (i) and (ii) do not contradict each other.
This shows (3)--(5).

\smallskip

In order to prove (6) we apply Theorem~\ref{thm:structureofgenselmer} with the set $Q$ constructed above. Indeed, elliptic curves are principally polarized and (4) implies that $E_{p^n}(k_v) \cong \bZ/p^n\bZ \times \bZ/p^n\bZ$.
\end{proof}

\section{Intersection with the maximal divisible subgroup}  \label{sec:maxdiv}

We are ultimately interested in understanding the intersection of $\Sha(A/k)$ with the maximal divisible subgroup of the Weil--Ch\^atelet group $\Div\big(\rH^1(k,A)\big)$. We proceed by analyzing one prime at a time. 


\subsection{Using the Euler characteristic formula}

Recall from \eqref{eq:pparth1localdiv} that $\Div\big(\rH^1(k,A)_{p^\infty}\big)$ is locally trivial at all primes $v$ which do not divide $p$. This implies that the image of $\rH^1_{\Sel_p}(k,A_{p^\infty})$ in 
$\rH^1(k,A)$ contains  $\Div\big(\rH^1(k,A)_{p^\infty}\big)$, and together with the analysis in 
Section \S\ref{sect:dualCond} gives us a way of exploring the intersection of $\Sha(A/k)_{p^\infty}$ with  $\Div\big(\rH^1(k,A)_{p^\infty}\big)$. We will follow this path first for a general abelian variety and then in the case of elliptic curves.

\begin{lem}\label{lem:trivp}
Let $A/k$ be an abelian variety, and let $p$ be a rational prime. If $\rH^1_{\Sel}(k,A^t_{p^n})$ equals $\rH^1_{\Sel^p}(k,A^t_{p^n})$ for almost all $n$ then $\Div\big(\rH^1(k,A)_{p^\infty}\big) \not= 0$
but  its intersection with $\Sha(A/k)_{p^\infty}$ agrees with $\Div\big(\Sha(A/k)_{p^\infty}\big)$.
\end{lem}

\begin{rmk}
The assumption $\rH^1_{\Sel}(k,A^t_{p^n}) = \rH^1_{\Sel^p}(k,A^t_{p^n})$ for almost all $n$ implies by \eqref{eq:selmerdef} that 
\[
A^t(k) \subseteq \bigcap_{n \geq 1} p^n A^t(k_v)
\]
where $v \mid p$ is a place of $k$. Therefore $A^t$ and, being isogenous, also $A$ have trivial algebraic rank.
\end{rmk}

\begin{proof}[Proof of Lemma~\ref{lem:trivp}] 
Consider the sequence \eqref{eq:selmerduality} for $M= A_{p^n} $, $\bL=\Sel_p$ and $\bL_0=\Sel$. Then, since $\bL^\ast=\Sel^p$, our assumption $\rH^1_{\Sel}(k,A^t_{p^n})=\rH^1_{\Sel^p}(k,A^t_{p^n})$ for $n \gg 0$ leads to the exact sequence
\begin{equation}\label{loc-surj}
0\rightarrow \rH^1_{\Sel}(k,A_{p^n})\rightarrow \rH^1_{\Sel_p}(k,A_{p^n}) \rightarrow \prod_{\fp|p} \rH^1(k_\fp, A)_{p^n}\rightarrow 0.
\end{equation}
Local Tate duality implies that we have the following group isomorphism
\begin{equation} \label{loc-size}
\prod_{\fp|p} \rH^1(k_\fp, A)_{p^n} \simeq (\bZ/p^n\bZ)^d \oplus \prod_{\fp|p} A(k_\fp)_{p^\infty}/{p^n}
\end{equation}
where $d=[k:\bQ] \cdot \dim(A)$ and the order of $\prod_{\fp|p} A(k_\fp)_{p^\infty}/{p^n}$ is independent of $n$ for $n\gg 0$.  In the limit for $n\to \infty$ the exact sequence (\ref{loc-surj}) becomes the exact sequence
\begin{equation}\label{eq:limitloc-surj}
0\rightarrow \rH^1_{\Sel}(k,A_{p^\infty})\rightarrow \rH^1_{\Sel_p}(k,A_{p^\infty}) \rightarrow \prod_{\fp|p} \rH^1(k_\fp, A)_{p^\infty}\rightarrow 0,
\end{equation}
and local Tate duality implies an isomorphism
\[
\prod_{\fp|p} \rH^1(k_\fp, A)_{p^\infty} = \prod_{\fp|p} \Hom\big(A(k_\fp),\bQ_p/\bZ_p\big) \cong (\bQ_p/\bZ_p)^d \oplus \prod_{\fp|p} A(k_\fp)_{p^\infty}
\]
Let $D_n \simeq (\bZ/p^n\bZ)^d$ be the intersection of the maximal divisible subgroup 
\[
D = \Div\big(\prod_{\fp|p} \rH^1(k_\fp, A)_{p^\infty}\big) \simeq  (\bQ_p/\bZ_p)^d
\]
with the $p^n$-torsion subgroup $\prod_{\fp|p} \rH^1(k_\fp, A)_{p^n}$.  Since (\ref{loc-surj}) is an exact sequence of finite length $\bZ/p^n\bZ$-modules, the set $S_n$ of partial splittings $s_n: D_n \to  \rH^1_{\Sel_p}(k,A_{p^\infty})$ is finite and non-empty. Restriction defines a map $S_{n+1} \to S_n$, and a well-known compactness argument shows that $\varprojlim_{n} S_n$ is non-empty. An element $s$ of the projective limit is nothing but a partial splitting $s: D \to  \rH^1_{\Sel_p}(k,A_{p^\infty})$ of (\ref{eq:limitloc-surj}). We conclude that the sequence
\begin{equation}\label{eq:maxdivlimitloc-surj}
0\rightarrow \Div\big(\rH^1_{\Sel}(k,A_{p^\infty})\big) \rightarrow \Div\big(\rH^1_{\Sel_p}(k,A_{p^\infty}) \big) \rightarrow \Div\big(\prod_{\fp|p} \rH^1(k_\fp, A)_{p^\infty}\big) \rightarrow 0.
\end{equation}
is split exact. Moreover, if we map to $\rH^1(k,A)_{p^\infty}$, in view of the discussion at the beginning of this section we find an exact sequence
\[
0 \to \Sha(A/k)_{p^\infty} \to \im\big(\rH^1_{\Sel_p}(k,A_{p^\infty}) \to \rH^1(k,A)_{p^\infty}\big) \to \prod_{\fp|p} \rH^1(k_\fp, A)_{p^\infty} \to 0,
\]
and again exploiting the partial splitting $s$ we find that
\[
0 \to \Div\big(\Sha(A/k)_{p^\infty}\big) \to \Div\big(\rH^1(k,A)_{p^\infty}\big) \to \Div\big(\prod_{\fp|p} \rH^1(k_\fp, A)_{p^\infty}\big) \to 0
\]
is also split exact. It follows that $\Div\big(\rH^1(k,A)_{p^\infty}\big)$ is nontrivial of corank at least 
\[
d = [k:\bQ]\cdot \dim(A)
\]
and that indeed the intersection of  $\Sha(A/k)_{p^\infty} \cap \Div\big(\rH^1(k,A)_{p^\infty}\big)$ lies in $\Div\big(\Sha(A/k)\big)$.
\end{proof}

\medskip

The \textbf{corank} of a $p$-primary torsion group $M = (\bQ_p/\bZ_p)^n \times M_0$ with finite $M_0$,  is the well-defined number $n = \coRk_{\bZ_p}(M)$.

\begin{prop} \label{prop:corankmaxdiv}
Let $A/k$ be an abelian variety with algebraic rank $r < d = \dim(A) \cdot [k:\bQ]$. Then $\Div(\rH^1(k,A))$ contains a copy of $(\bQ/\bZ)^{d-r}$.
\end{prop}
\begin{proof}
We can focus on the $p$-primary part for a prime number $p$ and have to compute the corank of $\Div(\rH^1(k,A)_{p^\infty})$. 
We set $s_p = \coRk_{\bZ_p}(\rH^1_{\Sel_p}(k,A_{p^\infty}))$,  $s = \coRk_{\bZ_p}(\rH^1_{\Sel}(k,A_{p^\infty}))$, and  $s^p = \coRk_{\bZ_p}(\rH^1_{\Sel^p}(k,A_{p^\infty}))$. These coranks are constant on isogeny classes, in particular they are the same for the dual $A^t$. 
Analyzing the asymptotic cardinality for $n \gg 0$ in \eqref{eq:selmerduality} with $M= A_{p^n} $, $\bL=\Sel_p$ and $\bL_0=\Sel$, we obtain 
\[
s_p - s^p = (s_p - s) + (s - s^p) = d.
\]
The exact sequence
\begin{equation} \label{eq:selmerfreeatpwithdivisiblecoefficents}
0 \to A(k) \otimes \bQ_p/\bZ_p \to \rH^1_{\Sel_p}(k,A_{p^\infty}) \to \ker \big (\rH^1(k,A)_{p^\infty}\to\prod_{v\nmid p}\rH^1(k_v,A)\big )  \to 0
\end{equation}
splits since $A(k) \otimes \bQ_p/\bZ_p$ is a divisible group. Therefore \eqref{eq:selmerfreeatpwithdivisiblecoefficents} remains exact upon applying  the functor $\Div(-)$.
Using the exact sequence 
\[
0 \to A(k) \otimes \bQ_p/\bZ_p \to \Div\big(\rH^1_{\Sel_p}(k,A_{p^\infty})\big) \to \Div\big(\rH^1(k,A)_{p^\infty}\big) \to 0
\]
we find
\begin{equation} \label{eq:corankfromula}
\coRk_{\bZ_p}\left(\Div\big(\rH^1(k,A)_{p^\infty}\big)\right) = s_p - r = d + s^p - r 
\end{equation}
which is $\geq d-r$ and proves the proposition.
\end{proof}


\subsection{Results for elliptic curves over \texorpdfstring{$\bQ$}{Q}} \label{sec:maxdivEllQ}

\begin{prop} \label{prop:trivialalgrank}
Let $E/\bQ$ be an elliptic curve of trivial algebraic rank. Then we have:
\begin{enumerate}
\item $\Div\big(\rH^1(\bQ,E)\big)$ contains a copy of $\bQ/\bZ$ and in particular is nontrivial.
\item If $\Sha(E/\bQ)$ is finite, then $\Div\big(\rH^1(\bQ,E)\big) \cong \bQ/\bZ$ and 
\[
\Sha(E/\bQ) + \Div\big(\rH^1(\bQ,E)\big) = \rH^1_{\divisor}(\bQ,E).
\]
\end{enumerate}
\end{prop}

\begin{proof}
(1) This is a special case of Proposition~\ref{prop:corankmaxdiv}.

(2) Fusing together the exact sequences \eqref{eq:pparth1localdiv} for all $p$ we obtain
\[
0 \to \Sha(E/\bQ) \to \rH^1_{\divisor}(\bQ,E) \to \bQ/\bZ \to 0
\]
which is exact since $\Sha(E/\bQ)$ was assumed finite and there is a copy of $\bQ/\bZ$ in  $\rH^1_{\divisor}(\bQ,E)$ by (1). The result follows at once.
\end{proof}

\begin{thm} \label{thm:trivintersection}
Let $E/\bQ$ be an elliptic curve of trivial analytic rank. Then we have:
\begin{enumerate}
\item[(1)] The intersection $\Sha(E/\bQ)_{p^\infty}$ with $\Div\big(\rH^1(\bQ,E)_{p^\infty}\big)$ is trivial for all odd primes $p$ of good reduction such that the $\Gal_\bQ$-representation on $E_p$ is irreducible.
\item[(2)]  $\Div\big(\rH^1(\bQ,E)\big) \cong \bQ/\bZ$ and the sum 
\[
\Sha(E/\bQ)_{p^\infty} \oplus \Div\big(\rH^1(\bQ,E)\big)_{p^\infty} = \rH^1_{\divisor}(\bQ,E)_{p^\infty}
\]
is direct for all $p$ as in (1).
\end{enumerate} 
\end{thm}

\begin{rmk}
We refer to the example in Section \S\ref{sec:selmerexample} for an elliptic curve $E/\bQ$ of trivial analytic rank with non-trivial $3$-torsion in $\Sha(E/\bQ) \cap \Div(\rH^1(\bQ,E))$.
\end{rmk}

\begin{proof}[Proof of Theorem~\ref{thm:trivintersection}] Since the analytic rank of $E/\bQ$ is trivial we know that $\Sha(E/\bQ)$ is finite (see \cite{Ko}), and hence (2) follows from (1) and Proposition~\ref{prop:trivialalgrank}.  We now prove (1) in several steps.

\smallskip

\textit{Step 1:}  We choose a quadratic imaginary extension $K=\bQ(\sqrt{-d})$ such that 
\begin{enumerate}
\renewcommand{\theenumi}{\roman{enumi}}
\renewcommand{\labelenumi}{(\theenumi)}
\item $K$ is distinct from the field of complex multiplication of $E$ in case $E$ has CM,
\item $K/\bQ$ is linearly disjoint from $\bQ(E_p)/\bQ$,
\item $d \geq 5$, 
\item $E/K$ has analytic rank $1$, 
\item all the primes dividing the conductor of $E/\bQ$ split,
\item $p$ is inert in $K/\bQ$. 
\end{enumerate}
Friedberg, and Hoffstein have shown that such an extension 
exists (see \cite{FH} Theorem B(2)).  

\smallskip

\textit{Step 2:} We fix $n\in \bN$ large enough so that
\begin{enumerate}
\renewcommand{\theenumi}{\roman{enumi}}
\renewcommand{\labelenumi}{(\theenumi)}
\item $p^{n-1}\Sha(E/K)_{p^\infty}=0$,
\item $p^n$ does not divide the basic Heegner point defined over $K$.
\end{enumerate}
This is possible since we know that $\Sha(E/K)$ is finite, and due to the nontriviality of this Heegner point, since the analytic rank of $E/K$ is $1$.

\smallskip

\textit{Step 3:} Since we have chosen $\bQ(E_p)$ and $K$ to be linearly disjoint over $\bQ$ the restriction yields an isomorphism
\[
G = \Gal(K(E_p)/K) \xrightarrow{\sim} \Gal(\bQ(E_p)/\bQ).
\]
Hence the assumption that $E_p$ is an irreducible  $\Gal_\bQ$-module, implies that it is also irreducible as a $\Gal_K$-module. Moreover, we have an isomorphism
\[
\Gal(K(E_p)/\bQ) = G \times \langle \sigma \rangle,
\]
where 
\[
\langle \sigma \rangle = \Gal(K(E_p)/\bQ(E_p)) \xrightarrow{\sim}  \Gal(K/\bQ).
\]

\smallskip

\textit{Step 4:} The image $\Gal(K(E_p)/K) \subseteq \GL(E_p)$ is not contained in a subgroup $S_3 \subseteq \GL(E_p)$, because $p=2$ was excluded, if $p=3$ then any $2$-dimensional $S_3$-representation is reducible, and for $p\geq 5$ we have
\[
   |\det(\Gal(K(E_p)/K))| =  |\det(\Gal(\bQ(E_p)/\bQ))| = [\bQ(\mu_p):\bQ] = p-1 > 2 = |\det(S_3)|.
\]

\textit{Step 5:} Let $L$ be the Selmer splitting field of $E$ over $K$ with respect to $p$, see Section~\ref{sec:Selmersplitting}. Because $E$ is defined over $\bQ$ and $L$ is characteristic for the base change of $E$ to $K$, we deduce that $L$ is Galois over $\bQ$. As in the 
proof of Proposition~\ref{prop:structuregenselmerellipticcurve} we conclude from Step 3 and 4 that $L$ and $K(E_{p^n})$ are linearly disjoint over $K(E_p)$. It follows from
\[
L \cap \bQ(E_{p^n}) = K(E_p) \cap \bQ(E_{p^n}) = \bQ(E_p)
\]
that the natural map to the fibre product
\begin{equation} \label{eq:galoisgroup}
\Gal(LK(E_{p^n})/\bQ) \xrightarrow{\sim} \Gal(L/\bQ) \times_{\Gal(\bQ(E_{p})/\bQ)} \Gal(\bQ(E_{p^n})/\bQ)
\end{equation}
is an isomorphism (here $LK(E_{p^n})$ denotes the compositum of $L$ and $K(E_{p^n})$).

\smallskip

\textit{Step 6:} Let $\tau \in \Gal(LK(E_{p^n})/\bQ) $ denote  a complex conjugation.  
The restriction of $\tau$ in 
\[
\Gal(K(E_p)/\bQ) =  G \times \langle \sigma \rangle \subseteq \GL(E_p) \times  \langle \sigma \rangle
\]
after a suitable choice of basis for $E_p$ reads 
\begin{equation} \label{eq:tau}
\matzz{1}{}{}{-1}  \cdot  \sigma.
\end{equation}

We have to understand $M=\Gal (L/K(E_p))$ as a $G$-module, more precisely as a $G \times \langle \sigma \rangle$-module. Recall from \eqref{eq:structureM} that naturally, and thus even as $G \times \langle \sigma \rangle$-modules, we have an inclusion
\[
M \inj \Hom(\rH^1_{\Sel}(K,E_p),E_p).
\]
 Here the product acts through projection to $\langle \sigma \rangle$ on $\rH^1_{\Sel}(K,E_p)$ and through projection to $G$ on $E_p$. Let $\chi_{K/\bQ}$ be the quadratic character associated to the extension $K/\bQ$. It follows ($p$ is odd) that as $G \times \langle \sigma \rangle$-module
\begin{equation} \label{eq:MintermsofEp}
M \cong (E_p)^a \oplus (E_p \otimes \chi_{K/\bQ})^b
\end{equation}
for some $a,b \in \bN$, because $E_p$ (and thus $E_p \otimes \chi_{K/\bQ}$) is irreducible and a submodule of a direct sum of simple modules is isomorphic to the direct sum over a subset of those summands.

\smallskip

Since $p$ is assumed to be odd, $M$ splits as a direct sum 
\[
M = M^+ \oplus M^-
\]
of its eigenspaces under $\tau$. We claim that $M^+$ generates $M$ as a $G$-module. 
By the decomposition \eqref{eq:MintermsofEp} it suffices to prove this claim for  the $G \times \langle \sigma \rangle$-modules $E_p$ and $E_p \otimes \chi_{K/\bQ}$. Due to \eqref{eq:tau} in both cases $\tau$ acts via a matrix conjugate to 
\[
\matzz{1}{}{}{-1}
\]
and the corressponding $+$-eigenspaces are nontrivial, so that these generate as a $G$-module by the irreducibility of $E_p$ (and thus $E_p \otimes \chi_{K/\bQ}$) as a $G$-module.

\smallskip

\textit{Step 7:}
Let $S\subseteq M=\Gal(L/K(E_p))$ be a generating set of $M$. 
We pick a finite set of rational primes $Q$ of $\bQ$ by choosing for every $h\in S$, using the Chebotarev density theorem,  a rational prime $\ell$ unramified in $LK(E_{p^n})/\bQ$ and coprime to the conductor of $E/\bQ$ such that 
\begin{enumerate}
\renewcommand{\theenumi}{\roman{enumi}}
\renewcommand{\labelenumi}{(\theenumi)}
\item $\Frob_\ell = \tau h \in \Gal(L/\bQ)$, and   
\item $\Frob_\ell = \tau \in  \Gal(\bQ(E_{p^n})/\bQ)$.
 \end{enumerate}
These requirements are free of contradictions  by the structure assertion \eqref{eq:galoisgroup}.

\smallskip

This set of auxiliary primes $Q$ satisfies the following properties:
\begin{enumerate}
\renewcommand{\theenumi}{\roman{enumi}}
\renewcommand{\labelenumi}{(\theenumi)}
\item $E$ has good reduction at $\ell\in Q$, 
\item $\ell \in Q$  is inert in $K/\bQ$,
\item $|E(K_\lambda)_{p^n}|=p^{2n}$ for every $\ell\in Q$ (with $\lambda$ the unique prime of $K$ above $\ell$ by (ii)),
\item the set of Frobenius elements $\Frob_{w} \in \Gal(L/K(E_p))$ where $w$ denotes a prime of $K(E_p)$ dividing $\ell$, when $\ell$ ranges over $Q$,  generates $\Gal(L/K(E_p))$ as a $\Gal(K(E_p)/K)$-module. 
\end{enumerate}
Indeed, properties (i)--(iii) are obvious and property (iv) holds by Step 6 because the $\tau$-eigenspace $M^+$  is generated by elements of the form 
\[
\Frob_w = (\Frob_\ell)^2 = (\tau h)^2 = \tau(h) \cdot h
\]
where $h \in S$. Here again, $p \not= 2$ is used.

\smallskip

We can then apply Theorem~\ref{thm:structureofgenselmer} with the auxiliary set $Q$ by now viewing $Q$ as a set of primes of $K$ (since there is a unique prime $\lambda$ of $K$ above each $\ell \in Q$). It follows that
\[
\rH^1_{\Sel_{Q}}(K,E_{p^n})\simeq (\bZ/p^n\bZ)^{2t}
\]
where $t=\#Q$. 

\smallskip

\textit{Step 8:} We now argue as in Theorem $1.1.7$ of  \cite{ciperianiwiles:solvpoints}. Since by Step 1 (ii) we still have $E_p(K) = 0$ and by Step 1 (iv) $E(K)$ contains a non-torsion Heegner point we have 
\[
E(K)/p^nE(K) \cong \bZ/p^n\bZ
\]
and the exact sequence
\[
0 \to E(K)/p^nE(K) \to \rH^1_{\Sel}(K,E_{p^n}) \to \Sha(E/K)_{p^n} \to 0
\]
splits. By the elementary divisor theorem for $\rH^1_{\Sel}(K,E_{p^n}) \subseteq \rH^1_{\Sel_{Q}}(K,E_{p^n})$  and Step 7, and by the large enough choice of $n$ in Step 2(i) we find $m_1 \leq  \ldots \leq m_{2t-1} < n$ such that 
\[
\rH^1_{\Sel}(K,E_{p^n}) \cong E(K)/p^nE(K)  \oplus \bZ/p^{m_1}\bZ \oplus \ldots \oplus \bZ/p^{m_{2t-1}}\bZ.
\]
Again by the choice of $n$ we know
\[
\rH^1_{\Sel}(K,E_{p^n})  \supset \bZ/p^{m_1}\bZ \oplus \ldots \oplus \bZ/p^{m_{2t-1}}\bZ \xrightarrow{\sim} \Sha(E/K)_{p^n} = \Sha(E/K)_{p^\infty},
\]
meaning that we can construct $\Sha(E/K)_{p^\infty}$ by constructing $2t-1$ independent ramified classes in 
$\rH^1_{\Sel}(K,E_{p^n}) \subseteq \rH^1_{\Sel_{Q}}(K,E_{p^n})$. 

\smallskip

\textit{Step 9:} 
Our choice of primes $\ell \in Q$ in Step 7 allows us to construct Kolyvagin classes  
 \[
 \{c_{\ell_1}, \ldots, c_{\ell_t}, c_{\ell_1\ell_2}, \ldots, c_{\ell_1\ell_t}\} \subseteq \rH^1_{\Sel_{Q}}(K,E_{p^n}),
 \]
 where $Q=\{\ell_1, \ldots \ell_t\}$. The assumption (ii) of Step 2 allows us to refine our choice of primes $Q$ in Step 7, as in section \S$1.4$ of  \cite{ciperianiwiles:solvpoints}, so that  
\begin{itemize}
\item $c_{\ell_1}$ is ramified at $\ell_1$, and 
\item $c_{\ell_i}$ and $c_{\ell_1\ell_i}$ are ramified at $\ell_i$ for every $i\geq 2$.
\end{itemize}
Finally, since  the analytic rank of $E/\bQ$ is trivial we know that the complex conjugation $\tau$ fixes $c_{\ell_i}$ and $\tau c_{\ell_1\ell_i}= -c_{\ell_1\ell_i}$.
This implies that we have constructed $2t-1$ independent ramified classes and that $\rH^1_{\Sel}(\bQ,E_{p^n})$ is contained in the subgroup of $\rH^1(K, E_{p^n})$ generated by $\{c_{\ell_1}, \ldots, c_{\ell_t}\}$. The assumption that $p$ is inert in $K/\bQ$ implies that $p$ splits completely in $K[\ell_i] /K$, where $K[\ell_i]$ denotes the ring class field of conductor $\ell_i$.  Consequently, the classes $\{c_{\ell_1}, \ldots, c_{\ell_t}\}$ are trivial at $p$ and hence
\[
\rH^1_{\Sel}(\bQ,E_{p^n})\hookrightarrow \rH^1_{\Sel^{p}}(K,E_{p^n}).
\]

Since $p$ is odd and $K/\bQ$ is a quadratic extension it follows that  
\[
\rH^1_{\Sel}(\bQ,E_{p^n})=\rH^1_{\Sel^{p}}(\bQ,E_{p^n}).
\]
Then, using Lemma~\ref{lem:trivp}, we see that the intersection $\Sha(E/\bQ)_{p^\infty}$ with $\Div\big(\rH^1(\bQ,E)_{p^\infty}\big)$ is a subgroup of $\Div\big(\Sha(E/\bQ)_{p^\infty}\big)$ which is trivial since $\Sha(E/\bQ)$ is finite.
\end{proof}

\begin{cor}
There are infinitely many elliptic curves $E/\bQ$ of trivial analytic rank such that the intersection of the Tate-Shafarevich group $\Sha(E/\bQ)$ and $\Div\big(\rH^1(\bQ,E)\big)$ is a $2 \rN$-power torsion group, where $\rN$ is the conductor of $E/\bQ$
\end{cor}

\begin{proof}
Recall that a Serre elliptic curve is an elliptic curve $E/\bQ$ such that the image of the product over all $p$ of the $p$-adic representations associated to the elliptic curve
\[
\rho_{E} : \Gal_\bQ \to \prod_p \GL(E_{p^\infty}) \cong \GL_2(\hat{\bZ})
\]
has index $2$ (the minimal possible index as noticed by Serre, see \cite{serre:pointsrationnels} Proposition~22).

Jones \cite{Jo10} and Zywina \cite{Zy10} consider the set of elliptic curves $E_{a,b}$ in the form 
\[
Y^2=X^3+aX+b
\]
such that $a$ and $b$ are integers of $\bQ$ that lie in the box defined by $h(a,b)< x$ (where $h$ denotes a height on such pairs). They show that the ratio of the cardinality of the subset of Serre elliptic curves 
by the total number of curves  in the box $h(a,b)< x$ approaches $1$ as $x$ goes to infinity. 

Clearly, for a Serre elliptic curve the $\Gal_\bQ$-representation $E_p$ is irreducible for all primes $p$. Hence in this sense  for 'most' elliptic curves $E/\bQ$, the $\Gal_\bQ$-representation $E_p$ is irreducible for every prime $p\geq 2$.

\smallskip

Friedberg and Hoffstein (see \cite{FH} Theorem B (1)) show that for every elliptic curve $E/\bQ$ there are infinitely many quadratic twists $E'/\bQ$ of trivial analytic rank. Observe that if $E_p$ is irreducible then so is $E'_p$. Hence we have infinitely many elliptic curves $E'$ of trivial analytic rank, irreducible $E'_p$ for every prime $p$. It then follows that the $p$-primary part of $\Sha(E/\bQ)$ and $\Div\big(\rH^1(\bQ,E)\big)$ intersect trivially for every odd prime $p$ of good reduction and the corollary follows.
\end{proof}

\begin{prop} \label{prop:rkalgpositive}
Let $E/\bQ$ be an elliptic curve of non-trivial algebraic rank. Then 
\[\Div\big(\rH^1(\bQ,E)\big) = 
\Div\big(\Sha(E/\bQ)\big).
\] 
\end{prop}

\begin{proof}
It suffices to argue for the $p$-primary part for every prime number $p$. Using \eqref{eq:selmerduality} for $M= E_{p^n} $, $\bL=\Sel_p$ and $\bL_0=\Sel$ together with local Tate duality, we see that 
\[
\frac{|\rH^1_{\Sel_p}(\bQ,E_{p^n})|}{|\rH^1_{\Sel}(\bQ,E_{p^n})|}  = |\rH^1(\bQ_p,E)_{p^n}|
\cdot \frac{|\rH^1_{\Sel^p}(\bQ,E_{p^n})|}{|\rH^1_{\Sel}(\bQ,E_{p^n})|}
\]
\[
 \leq 
\frac{|E(\bQ_p)/p^nE(\bQ_p)|}{|\im\big(E(\bQ)/p^nE(\bQ) \to E(\bQ_p)/p^nE(\bQ_p)\big)|}
= |E(\bQ_p)/(p^nE(\bQ_p) + E(\bQ))|
\]
which is bounded independently of $n$ since the algebraic rank of $E/\bQ$ is non-trivial. In view of the discussion at the beginning of this section it then follows that
 $\Div\big(\rH^1(\bQ,E)_{p^\infty}\big)\subseteq \Div\big(\Sha(E/\bQ)\big)_{p^\infty}$. The other inclusion is clear.
\end{proof}

\medskip

Proposition~\ref{prop:rkalgpositive} above can also be deduced from \cite{bashmakov:survey} Theorem 7 and can essentially be found in \cite{harariszamuely:galsec} Corollary 4.2. It has the following consequence.

\begin{cor} \label{cor:analyticrank1}
Let $E/\bQ$ be an elliptic curve of analytic rank $1$. Then $\Div\big(\rH^1(\bQ,E)\big)$ is trivial.
\end{cor}
\begin{proof}
Kolyvagin \cite{Ko} has shown that $\Sha(E/\bQ)$ is finite for all elliptic curves $E/\bQ$ of analytic rank $1$ and that the algebraic rank is equal to $1$ as well. Hence by Proposition~\ref{prop:rkalgpositive} we find $\Div\big(\rH^1(\bQ,E)_{p^\infty}\big)\subseteq \Div\big(\Sha(E/k)\big)_{p^\infty}=0$ for all primes $p$.
\end{proof}


\section{An example: the Jacobian of the Selmer curve} \label{sec:selmerexample}
We are now discussing in detail the plane cubic 
\[
S = \{3X^3 + 4Y^3 + 5Z^3 = 0\}
\]
describing Selmer's curve of genus $1$ that violates the Hasse principle. Its Jacobian $E=\Pic_S^0$ is an elliptic curve over $\bQ$ of analytic rank $0$ 
given by the homogeneous equation
\begin{equation} \label{eq:selmercurve}
X^3 + Y^3 + 60Z^3 = 0
\end{equation}
with $[1:-1:0]$ as its origin.
The curve $E$ has  Mordell-Weil group $E(\bQ)=0$, see \cite{cassels:lectures}  \S18 Lemma 2, and $3$-torsion in an exact sequence
\begin{equation} \label{eq:3torsionofselmerexample}
0 \to \mu_3 \to E_3 \to \bZ/3\bZ \to 0,
\end{equation}
which splits over a field $k/\bQ$ if and only if $60$ is a cube in $k$. Note that \eqref{eq:3torsionofselmerexample} shows that with respect to the prime $p=3$ we are in the $S_3$-case excluded in Proposition~\ref{prop:structuregenselmerellipticcurve} (on top of $E_3$ not being irreducible).

\smallskip

The curve $S$, as a principal homogeneous space under $E$ describes a nontrivial $3$-torsion element of $\Sha(E/\bQ)$, see \cite{mazur:localtoglobal} I \S4 and \S9. 
Mazur and Rubin determine 
\[
\Sha(E/\bQ) \cong \bZ/3\bZ \times \bZ/3\bZ
\]
unconditionally,  see \cite{mazur:localtoglobal} Theorem 1 and \S9, and give a list of the isomorphism classes as genus $1$ curves of the principal homogeneous spaces\footnote{The element $W \in \Sha(E/\bQ)$ is isomorphic to $-W$ as a curve of genus $1$.} under $E$ representing nontrivial elements of $\Sha(E/\bQ)$ as, besides $S$, 
\[
S' = \{X^3 + 5Y^3 + 12Z^3 = 0\},
\]
\begin{equation} \label{eq:shaselmer}
S'' = \{X^3 + 4Y^3 + 15Z^3 = 0\},
\end{equation}
\[
S''' = \{X^3 + 3Y^3 + 20Z^3 = 0\}.
\]
Kummer theory and the vanishing $E(\bQ) = 0$ yield isomorphisms
\[
\rH^1_{\Sel}(\bQ,E_3) \xrightarrow{\sim}  \rH^1_{\Sel}(\bQ,E_{3^n})  \xrightarrow{\sim} \Sha(E/\bQ).
\]
The $E_3$-torsor associated to an element of $\Sha(E/\bQ)$ is given by the zero set of $XYZ=0$ in the description given by \eqref{eq:selmercurve} and \eqref{eq:shaselmer}. The Selmer group trivial at $3$ can now be determined via the defining exact sequence
\[
0 \to \rH^1_{\Sel^3}(\bQ,E_3) \to \rH^1_{\Sel}(\bQ,E_3) \to \rH^1(\bQ_3,E_3)
\]
since we only need to check whether $XYZ=0$ on $S,S',S''$ or $S'''$ has a $\bQ_3$-point. This is true for $S$ and false for $S',S'',S'''$ as can be deduced from $\bQ_3^\ast/(\bQ_3^\ast)^3 = \langle 2,3\rangle$ with $10$ being a cube in $\bQ_3^\ast$. We conclude that
\[
\Sha(E/\bQ) \supset \langle [S] \rangle = \rH^1_{\Sel^3}(\bQ,E_3) = \rH^1_{\Sel^3}(\bQ,E_{3^n}) \cong \bZ/3\bZ.
\]
Since $E(\bQ_3)$ has no $3$-torsion, we find by local Tate duality 
\[
\rH^1(\bQ_3,E)_{3^\infty}  = \Hom(E(\bQ_3), \bQ_3/\bZ_3) = \bQ_3/\bZ_3,
\]
so that counting in \eqref{eq:selmerduality} for $M= E_{3^n} $, $\bL=\Sel_3$ and $\bL_0=\Sel$, leads to
\begin{equation} \label{eq:sizeofselmerfreeat3}
|\rH^1_{\Sel_3}(\bQ,E_{3^n})| = 3^{n+1}.
\end{equation}
The map $\rH^1_{\Sel_3}(\bQ,E_{3^n}) \to \rH^1(\bQ,E)$ is injective with image the $3^n$-torsion of the subgroup
\[
\rH^1_{\divisor}(\bQ,E)_{3^\infty} 
= \ker\big(\rH^1(\bQ,E) \to \prod_{v \not= 3} \rH^1(\bQ_v,E)\big),
\]
and hence 
\[
\rH^1_{\Sel_3}(\bQ,E_{3^\infty}) = \rH^1_{\divisor}(\bQ,E)_{3^\infty} \cong \bZ/3\bZ \times \bQ_3/\bZ_3.
\]
Here we have used the elementary fact that for a finitely cogenerated torsion group $M$ the knowledge of the orders of $n$-torsion $M_n$ for all $n$ uniquely determines the structure as an abstract group: the sizes are given by \eqref{eq:sizeofselmerfreeat3}, while 
\[
\rH^1_{\Sel_3}(\bQ,E_{3^n}) = \left(\rH^1_{\Sel_3}(\bQ,E_{3^\infty})\right)_{3^n}
\]
shows that we have indeed counted the exact $3^n$-torsion of a finitely cogenerated torsion group.

It follows that $\Sha(E/\bQ)$ meets the maximal divisible group of $\rH^1(\bQ,E)$ nontrivially, namely
\[
\Sha(E/\bQ) \cap \Div(\rH^1(\bQ,E)) = \Sha(E/\bQ) \cap \Div(\rH^1_{\divisor}(\bQ,E)) \cong \bZ/3\bZ,
\]
which shows the limitation of Theorem~\ref{thm:trivintersection} towards $p=3$: the prime $p=3$ happens to be a prime of additive bad reduction for $E/\bQ$ and $E_3$ is a reducible $\Gal_\bQ$-module.

\smallskip

We proceed to determine, which of  the principal homogeneous spaces $S,S',S''$ or $S'''$ generates the intersection $\Sha(E/\bQ) \cap \Div(\rH^1(\bQ,E))$, thereby preparing the answer to  \cite{stix:bm} Question~49 given in Remark~\ref{rmk:selmercurvesplits} below and completing the proof of Theorem B of the introduction.

\begin{prop} \label{prop:selmercurvedivisible}
The intersection $\Sha(E/\bQ) \cap \Div(\rH^1(\bQ,E))$ is generated by the class of the Selmer curve $S$.
\end{prop}
\begin{proof}
The Jacobian $E$ of the Selmer curve has good reduction outside $\{2,3,5\}$. 
We set 
\[
U = \Spec(\bZ[1/30])
\]
and consider $E$ by abuse of notation as the relative elliptic curve $E/U$. 
Since for $v = 2,5$ neither a cube root of $60$ nor $\zeta_3$ is contained in $\bQ_v$, we have $E_3(\bQ_v) =0$ and thus by local Tate-duality 
\[
\rH^1(\bQ_v,E_3) = \rH^1(\bQ_v,E_9) = 0
\]
for $v = 2,5$. We conclude by Corollary~\ref{cor:comparegeneralizedselmer} and 
\eqref{eq:compactlocalisation}
that 
\begin{eqnarray*}
\rH^1(U,E_3)  = & \rH^1_{\Sel_3}(\bQ,E_3) & = \Sha(E/\bQ) \cong \bZ/3\bZ \times \bZ/3\bZ, \\
\rH^1_{\rc}(U,E_3) = & \rH^1_{\Sel^3}(\bQ,E_3) & =  \langle [S] \rangle \cong \bZ/3\bZ, \\
\rH^1(U,E_9) = & \rH^1_{\Sel_3}(\bQ,E_9) & \cong \bZ/3\bZ \times \bZ/9\bZ.
\end{eqnarray*}

In order to decide which element of $\Sha(E/\bQ)$ generates the intersection with $\Div(\rH^1(\bQ,E))$, it suffices to check which classes are divisible by $3$ in $\rH^1(U,E_{9})$ which is the $9$-torsion of 
\[
\rH^1(U,E_{3^\infty})  = \rH^1_{\Sel_3}(\bQ,E_{3^\infty}) = \bZ/3\bZ \times \bQ_3/\bZ_3.
\]
This is controlled by a Bockstein map as we will explain now. The long exact cohomology sequence for 
\begin{equation} \label{eq:393Kummer}
0 \to E_3 \to E_9 \to E_3 \to 0
\end{equation}
reads
\[
0 \to \rH^1(U,E_3) \to \rH^1(U,E_9) \xrightarrow{"3 \cdot"} \rH^1(U,E_3) \xrightarrow{\beta} \rH^2(U,E_3) \to 0
\]
where the \textbf{Bockstein} map $\beta:\rH^1(U,E_3) \to \rH^2(U,E_3)$ is surjective by counting and Artin--Verdier duality 
\[
\rH^2(U,E_3) = \Hom(\rH^1_{\rc}(U,E_3),\bQ/\bZ) \cong \bZ/3\bZ
\]
induced by the Weil pairing $e : E_3 \otimes E_3 \to \mu_3$, see  \cite{zink:appendix} Theorem 3.2 
and \cite{zink:appendix} Lemma 3.2.2. 
We introduce the \textbf{Weil--Bockstein} pairing
\[
\langle-,-\rangle_{\rm WB} \ : \ \rH^1(U,E_3) \times \rH^1_{\rc}(U,E_3) \to \bZ/3\bZ
\]
\[
\langle a,b \rangle_{\rm WB} = e_\ast\big(\beta(a) \cup b\big).
\]
With the  "forget support map" 
\[
f^i : \rH^i_{\rc}(U,E_3) \to \rH^i(U,E_3)
\]
and the class of the Selmer curve $s \in \rH^1_{\rc}(U,E_3)$ the claim of the proposition is equivalent to $\beta(f^1(s)) =0$ which boils down to the vanishing of 
\[
\langle f^1(s),s\rangle_{\rm WB}.
\]
Applying compactly supported cohomology to \eqref{eq:393Kummer} we find the 
\textbf{compact Bockstein}
\[
\beta_{\rc} : \rH^1_{\rc}(U,E_3) \to \rH^2_{\rc}(U,E_3)
\]
which is adjoint to $\beta$ with respect to Artin--Verdier duality sponsored by the Weil pairing. Moreover, the forget support maps $f^1$ and $f^2$ are also adjoints. We can thus compute
\[
\langle f^1(s),s\rangle_{\rm WB} = e_\ast \big( f^1(s)) \cup \beta_{\rc}(s)\big) 
= e_\ast \big(s \cup f^2(\beta_{\rc}(s))\big) = e_\ast \big(s \cup \beta(f^1(s))\big) = - \langle f^1(s),s \rangle_{\rm WB}
\]
due to anti-symmetry of the Weil-pairing, so that indeed $\langle f^1(s),s\rangle_{\rm WB} = 0$.
\end{proof}

\begin{rmk}
The proof of Proposition~\ref{prop:selmercurvedivisible} given above can be used to identify the pairing 
$\langle-,-\rangle_{\rm WB}$
with the restriction of the Cassels--Tate pairing on $\Sha(E/\bQ)$ in the form of the "Weil--pairing definition" of Poonen and Stoll \cite{PS} \S12.2. 
\end{rmk}

\begin{rmk} \label{rmk:selmercurvesplits}
Proposition~\ref{prop:selmercurvedivisible} answers \cite{stix:bm}  Question~49. The Selmer curve $S/\bQ$
provides an explicit example for a smooth projective curve with no rational point over $\bQ$, but rational points everywhere locally, and nevertheless split fundamental group extension 
\[
1 \to \pi_1(\ov{S}, \bar s) \to \pi_1(S,\bar s) \to \Gal_k \to 1
\]
due to \cite{stix:habil}  Corollary 177. The splitting obviously does not come from a rational point, see \cite{stix:habil} \S13 for the 
context of the section conjecture of anabelian geometry.
\end{rmk}

\begin{rmk} We will show in a subsequent paper \cite{ciperianistix:cassels}, when we discuss the divisibility question in the sense of Cassels, that $\Sha(E/\bQ)$, and consequently $\rH^1_{\divisor}(\bQ,E)_{3^\infty}$, is contained in $\divisor(\rH^1(\bQ,E))$.  Therefore the Jacobian of the Selmer curve gives a concrete example where the filtration \eqref{eq:filtration} of the introduction has all its steps nontrivial.
\end{rmk}



\begin{thebibliography}{NSW08}

\bibitem[AI82]{andersonihara:proell} 
Anderson, G., Ihara, Y., 
Pro-$\ell$ branched coverings of $\bP^1$ and higher circular $\ell$-units, 
\textit{Ann.\ of Math.} (2) \textbf{128} (1988), no.\ 2, 271--293. 

\bibitem[Ba64]{bashmakov:div} 
Bashmakov, M. I., 
On the divisibility of principal homogeneous spaces over Abelian varieties, 
\textit{Izv.\ Akad.\ Nauk SSSR Ser.\ Mat.\ } \textbf{28} (1964), 661--664. 

\bibitem[Ba72]{bashmakov:survey} 
Bashmakov, M. I., 
The cohomology of abelian varieties over a number field, (English translation) 
\textit{Russian Math. Surveys} \textbf{27} (1972), no.\ 6, 25--70.

\bibitem[Ca62]{cassels:III} 
Cassels, J. W. S.,
Arithmetic on curves of genus $1$. III. The Tate-Shafarevic and Selmer groups,
\textit{Proc. London Math. Soc. }(3) \textbf{12}  1962, 259--296.

\bibitem[Ca91]{cassels:lectures} 
Cassels, J. W. S., 
\textit{Lectures on elliptic curves},  
London Mathematical Society Student Texts  \textbf{24}, Cambridge University Press, 1991, vi+137pp. 

\bibitem[\c{C}S12]{ciperianistix:cassels}
\c{C}iperiani, M., Stix, J., 
Weil--Ch\^atelet divisible elements in Tate--Shafarevich groups II: On a question of Cassels,
preprint, 2012.

\bibitem[\c{C}W08]{ciperianiwiles:solvpoints}  
\c{C}iperiani, M., Wiles, A., 
Solvable points on genus one curves, 
\textit{Duke Mathematical Journal} \textbf{142} (2008), no.\ 3, 381--464.

\bibitem[FH95]{FH} 
Friedberg, S., Hoffstein, J., 
Nonvanishing theorems for automorphic $L$-functions on ${\rm GL}(2)$, 
\textit{ Ann. of Math.} (2)   \textbf{142}  (1995),  no.\ 2, 385--423.
   
\bibitem[HS09]{harariszamuely:galsec} 
Harari, D., Szamuely, T., 
Galois sections for abelianized fundamental groups, with an Appendix by E.~V. Flynn, 
\textit{Math.\ Ann.} \textbf{344} (2009), no.\ 4, 779--800. 

\bibitem[Ja88]{jannsen:continuous} 
Jannsen, U., 
Continuous \'etale cohomology, {\it Math.\ Annalen} {\bf 280} (1988), 207--245.

\bibitem[Jo10]{Jo10} 
Jones, N., 
Almost all elliptic curves are Serre curves, 
\textit{Trans. Amer. Math. Soc.} \textbf{362}  (2010),  no.\ 3, 1547--1570.

\bibitem[Ko90]{Ko} 
Kolyvagin, V. A., 
Euler systems, 
\textit{The Grothendieck Festschrift}, Vol.\ II,  Progr. Math., 87, Birkhauser Boston, Boston, MA, 1990, 435--483.

\bibitem[Ma78]{mazur:eisenstein} 
Mazur, B., 
Modular curves and the Eisenstein ideal, \textit{IHES Publ.\ Math.} \textbf{47}, (1978), 33--186.

\bibitem[Ma93]{mazur:localtoglobal} 
Mazur, B., 
On the passage from local to global in number theory,
\textit{Bull.\ Amer.\ Math.\ Soc.} \textbf{29} (1993), no.\ 1, 14--50. 

\bibitem[Mi86]{Milne} 
Milne, J.S., 
\textit{Arithmetic duality theorems}, 
Perspectives in Mathematics, 1. Academic Press, Inc., Boston, MA, 1986.

\bibitem[Ne75]{neumann:pclosed} 
Neumann, O., 
$p$-closed algebraic number fields with bounded ramification, 
\textit{Izv.\ Akad.\ Nauk SSSR Ser.\ Mat.\ } \textbf{39} (1975), no.\ 2, 259--271, 471.

\bibitem[NSW08]{nsw}  
Neukirch, J., Schmidt, A., Wingberg, K., 
\textit{Cohomology of number fields},
second edition, Grundlehren der Mathematischen Wissenschaften \textbf{323}, Springer, 2008, xvi+825pp.

\bibitem[PS99]{PS} 
Poonen, B., Stoll, M., 
The {C}assels-{T}ate pairing on polarized abelian varieties, 
\textit{Ann.\ of Math.\ (2)} \textbf{150} (1999), no.\ 3, 1109--1149.

\bibitem[Ri97]{ribet:sstgalois} 
Ribet, K., 
Images of semistable Galois representations, 
\textit{Pacific J. of Math.} \textbf{181} (1997), no.\ 3, 277--297.

\bibitem[Se72]{Serre} 
Serre, J.-P., 
Propri\'et\'es galoisiennes des points d'ordre fini des courbes elliptiques, 
\textit{Invent. Math.}  \textbf{15}  (1972), no.\ 4, 259--331.  

\bibitem[Se79]{serre:pointsrationnels} 
Serre, J.-P., 
Points rationnels des courbes modulaires $X_0(N)$ [d'apr\`es Barry Mazur], 
\textit{S\'eminaire Bourbaki} (1977/78), expos\'e \textbf{511}, Lecture Notes in Math. \textbf{710} (1979), 89--100. 

\bibitem[Sh63]{shafarevich} 
Shafarevich, I.~R., 
Extensions with given ramification points (russian),
\textit{Publ.\ Math.\ IHES} \textbf{ 18} (1963), 71--95, translated in \textit{Amer.\ Math.\ Soc.\ Transl.} \textbf{59} (1966), 128--149.
  
\bibitem[Sx11]{stix:bm} 
Stix, J., 
The Brauer--Manin obstruction for sections of the fundamental group, 
\textit{Journal of Pure and Applied Algebra} \textbf{215} (2011), no.\ 6, 1371--1397.

\bibitem[Sx12]{stix:habil} 
Stix, J., 
\textit{Rational Points and Arithmetic of Fundamental Groups: Evidence for the Section Conjecture}, 
Springer Lecture Notes in Mathematics \textbf{2054}, Springer, 2012, xx+247 pp.

\bibitem[Zi78]{zink:appendix} 
Zink, Th., 
\'Etale cohomology and duality in number fields,
Appendix 2 in:   Haberland, K., {\it Galois Cohomology of Algebraic Number Fields}, VEB  Deutscher Verlag der Wissenschaften, 1978.

\bibitem[Zy10]{Zy10}
Zywina, D., 
Elliptic curves with maximal Galois action on their torsion points, 
\textit{ Bull.\ Lond.\ Math.\ Soc.} \textbf{ 42}  (2010),  no.\ 5, 811--826.


\end{thebibliography}
\end{document}